\documentclass[11pt,reqno,oneside]{amsart}
\usepackage[english]{babel}
\usepackage[utf8]{inputenc}
\usepackage[mathscr]{eucal}
\usepackage{xcolor}
\usepackage{amsfonts}   
\usepackage{amssymb}
\usepackage{latexsym}
\usepackage{enumerate}
\usepackage[colorlinks,citecolor=purple,linkcolor=blue]{hyperref}

\usepackage{geometry}

\let\nc\newcommand

\nc{\la}{\label}

\newtheorem{theorem}{Theorem}[section]
\newtheorem{corollary}[theorem]{Corollary}
\newtheorem{lemma}[theorem]{Lemma}
\newtheorem{proposition}[theorem]{Proposition}
\theoremstyle{definition}
\newtheorem{definition}[theorem]{Definition}
\newtheorem{example}[theorem]{Example}

\theoremstyle{remark}
\newtheorem{remark}[theorem]{Remark}

\def\k{\mathsf k}
\def\C{\mathbb C}

\newcommand{\Frac}{{\rm{Frac}}}

\newcommand{\Spec}{{\rm{Spec}}}

\newcommand{\Autk}{{\rm{Aut}}}

\newcommand{\Der}{{\rm{Der}}}

\newcommand{\into}{\,\,\hookrightarrow\,\,}

\def\k{\mathsf k}
\newcommand{\Tdeg}{{\rm Tdeg}}
\newcommand{\LD}{{\rm LD}}

\newcommand{\gr}{{\rm gr}}
\newcommand{\GK}{{\rm GK}}
\def\K{\mathcal K}
\def\KK{\mathfrak K}
\def\L{\mathcal L}
\def\M{\mathcal M}
\def\Q{\mathcal Q}

\def\cl{\operatorname{cl.Krull}}
\def\p{\mathfrak p}

\newcommand{\End}{\operatorname{End}}

\makeatletter
\newcommand{\colim@}[2]{%
  \vtop{\m@th\ialign{##\cr
    \hfil$#1\operator@font colim$\hfil\cr
    \noalign{\nointerlineskip\kern1.5\ex@}#2\cr
    \noalign{\nointerlineskip\kern-\ex@}\cr}}%
}
\newcommand{\colim}{%
  \mathop{\mathpalette\colim@{\rightarrowfill@\scriptscriptstyle}}\nmlimits@
}
\renewcommand{\varprojlim}{%
  \mathop{\mathpalette\varlim@{\leftarrowfill@\scriptscriptstyle}}\nmlimits@
}
\renewcommand{\varinjlim}{%
  \mathop{\mathpalette\varlim@{\rightarrowfill@\scriptscriptstyle}}\nmlimits@
}
\makeatother

\begin{document}

\title[Galois Rings: Ring-Theoretic Properties and Applications]{Galois Rings: Ring-Theoretic Properties and Applications to Coulomb Branches\\and Affine Hecke Algebras}

\author{Vyacheslav Futorny} \address{Shenzhen International Center for Mathematics, Southern University of Science and Technology, Shenzhen, China} 
  \email{vfutorny@gmail.com}
  \thanks{ V.F. is partially supported by the NSF of China grants 12350710178 and 12350710787.}
\author{Jonas T. Hartwig}
\address{Department of Mathematics, Iowa State University, Ames IA 50011, USA}
\email{jth@iastate.edu}
\urladdr{http://jthartwig.net}
\thanks{J.T.H. is partially supported by the United States Army Research Office grant W911NF-24-1-0058.}

\author{Erich C. Jauch}
\address{Department of Mathematics \& Physics, Westminster College (Missouri), Fulton MO 65251}
\email{erich.jauch@westminster-mo.edu}
\urladdr{http://ecjauch.com}

\author{Jo\~ao Schwarz}
\address{Shenzhen International Center for Mathematics, Southern University of Science and Technology, Shenzhen, China}
\email{jfschwarz.0791@gmail.com}


\begin{abstract}
Galois rings and Galois orders, introduced by Futorny and Ovsienko, are realized as subrings of fixed subrings of skew group (or monoid) rings and have numerous applications in the structure and representation theory of associative algebras. This paper consists of two parts. The first part investigates ring-theoretic properties that follow from the Galois ring structure alone. In particular, we establish natural conditions under which Galois rings are Ore domains or (semi)prime Goldie rings. We also study several ring-theoretic dimensions and combine the theories of Galois rings and PI-algebras to obtain new structural results.

In the second part, we apply these results, together with general techniques from ring theory, to affine Hecke algebras in the sense of Ginzburg, Kapranov, and Vasserot, as well as to spherical Coulomb branch algebras. In particular, we prove that these algebras are Jacobson semiprimitive, satisfy the Nullstellensatz, and determine several of their ring-theoretic dimensions. For affine Hecke algebras, we further prove that they satisfy the maximal Nullstellensatz, are integral over their centers, and determine their T-ideals of polynomial identities, PI-degrees, and PI-exponents. For spherical Coulomb branch algebras, we compute the Krull dimension, establish that they satisfy the Gelfand–Kirillov conjecture, and prove that they are not PI-algebras.
\end{abstract}

\maketitle

\tableofcontents

\section{Introduction}
The concepts of Galois rings and Galois orders were introduced by Futorny and Ovsienko in \cite{FO,FO2} to provide a suitable framework for the representation theory of certain infinite-dimensional Noetherian algebras. This theory unified the Gelfand--Tsetlin theory for $\mathfrak{gl}_n$, developed in \cite{DFO0,DFO}, with the representation theory of generalized Weyl algebras introduced by V. Bavula in \cite{Bavula}. It can also be viewed as a refinement of the general theory of Harish--Chandra modules, initiated in \cite{DFO} and further developed in \cite{Webster}, \cite{Fillmore}, \cite{FillHart}, and \cite{Schwarz2}. The underlying idea goes back to Block \cite{Block}: to study irreducible modules by exploiting a suitable embedding of the algebra into a skew group ring.

The theory of Galois orders has been successfully applied to the study of representations of generalized Weyl algebras \cite{BO,BavulaK}, finite $W$-algebras of type $A$ \cite{FMO}, invariants of certain rings of differential operators \cite{FS0,FS3,Schwarz2}, invariants of quantum groups \cite{FS2}, the alternating analogue of $U(\mathfrak{gl}_n)$ \cite{Jauch}, and OGZ algebras, their $q$-analogues, and parabolic versions of $U_q(\mathfrak{gl}_n)$ \cite{Hartwig}. In particular, the notions of principal and rational Galois orders were introduced in the latter paper. Moreover, important variations of these concepts, namely flag orders, were introduced in \cite{Webster}, where it was shown that spherical Coulomb branch algebras, defined in \cite{BFN} (see also \cite{lots}), are principal Galois orders. The Galois order realization of spherical Coulomb branch algebras was subsequently used in \cite{LW} to study spherical subalgebras of rational Cherednik algebras, and in \cite{lots} in a much more general setting. Further developments of the theory can be found in \cite{FGRZ}, \cite{MV}, \cite{Jauch2}, \cite{Hartwig2}, \cite{Hartwig3}, \cite{Fillmore}, \cite{Schwarz3}, and \cite{Schwarz2}.

Ring-theoretic aspects have always played an important role in the development of this theory. Galois orders generalize the classical theory of orders (see, e.g., \cite[Chapters 3 and 5]{McConnell}), where the denominator set is no longer required to be central.

The theory of Galois algebras has also provided a powerful tool for establishing the validity of the Gelfand--Kirillov conjecture \cite{Gelfand} (abbreviated here as GKC) and its $q$-analogue (see \cite[I.2.11, II.10.4]{Brown}, abbreviated q-GKC) for many classes of algebras; see, for example, \cite{FMO,FH,EFOS,H0,FS2,FS3,Jauch,Hartwig3,Schwarz2}. It also provides effective methods for studying the Gelfand--Kirillov dimension, the center, and maximal commutative subalgebras of algebras that admit realizations as Galois orders \cite{FO,Hartwig}. In particular, it was proved in \cite{Hartwig} that the Gelfand--Tsetlin subalgebra of $U_q(\mathfrak{gl}_n)$ is maximal commutative whenever $q$ is not a root of unity, thereby confirming a long-standing conjecture of Mazorchuk and Turowska \cite{MT}. More recently, the theory has also been used to study the $\mathbb{Z}$-graded polynomial identities of certain generalized Weyl algebras \cite{FKS}.

This paper consists of two parts. In the first part, we investigate ring-theoretic properties of Galois rings that follow solely from the Galois ring structure. In the second part, we apply these results, together with more general techniques from ring theory, to the (double) affine Hecke algebras constructed in \cite{GKV} and to spherical Coulomb branch algebras introduced in \cite{BFN}. Although spherical Coulomb branch algebras have been extensively studied in connection with symplectic singularities, symplectic duality, and representation theory \cite{Joel,Webster,lots}, their ring-theoretic properties have received comparatively little attention.

Throughout the paper, we work in the framework developed in \cite{Hartwig}, since it does not require the rings under consideration to be algebras over a base field. When the rings are assumed to be algebras over an algebraically closed field of characteristic zero, the frameworks of \cite{FO} and \cite{Hartwig} are essentially equivalent by \cite[Theorem~4.1]{Schwarz2}.

The paper is organized as follows. In Section~2, we recall the fundamentals of the theory of Galois rings from \cite{FO} and \cite{Hartwig} and extend several elementary results from \cite{FO} and \cite{Schwarz2} to the framework of \cite{Hartwig}. In particular, we discuss generalized Weyl algebras, maximal commutative subalgebras, and the center of a Galois ring. Section~3 forms the core of the paper. We show that, given a realization of an associative ring $U$ as a $\Gamma$-ring in a fixed subring $(\L * \M)^G$ of a skew monoid ring, several ring-theoretic properties of $U$ can be derived either from those of $(\L * \M)^G$ (see Theorems~\ref{localization4} and \ref{Ore}) or solely from the monoid $\M$ (for example, being an Ore domain, Theorem~\ref{U-is-Ore}, or a prime Goldie ring, Theorem~\ref{U-is-prime-Goldie}).

In Section~4, we combine the theories of Galois rings and PI-algebras. We prove that $U$ is a PI-algebra if and only if $\K$ is a PI-algebra, in which case they satisfy the same polynomial identities (Theorem~\ref{PI-0}). When $U$ is both a prime PI-algebra and a Galois ring, stronger structural results are obtained in Theorems~\ref{PI} and \ref{integrality}. In Theorem~\ref{usually-not-PI}, we develop a simple criterion for proving that a Galois ring is not a PI-algebra. As a first application, we show that the alternating analogue of $U(\mathfrak{gl}_n)$ introduced in \cite{Jauch} is not a PI-algebra. 
Finally, we apply the methods developed for Galois PI-algebras to the nilHecke algebra \cite{G18}. 

Section~5 is devoted to the study of five ring-theoretic dimensions of Galois rings: the classical Krull dimension $\cl$, the Gabriel--Rentschler Krull dimension $\KK$, the Gelfand--Kirillov dimension $\GK$, the Gelfand--Kirillov transcendence degree $\Tdeg$, and the lower transcendence degree $\LD$. We begin by proving a result of independent interest on the lower transcendence degree (Proposition~\ref{new-LD}). This is immediately applied in Corollary~\ref{useful-degrees} to the algebras $S_m^n(\k)$ and $Q_m^n(\k)$ (see Definition~\ref{main-chars}), which play a fundamental role in the theory of Galois rings (see, for example, \cite{Hartwig,Schwarz2}). We then prove Theorem~\ref{loc-alg}, which simplifies \cite[Theorem~6.1]{FO}. As further applications, we show that the noncommutative deformations of Kleinian singularities of type $D$ and the quantum enveloping algebra $U_q(\mathfrak{gl}_n)$, where $q$ is not a root of unity, are $\LD$-stable, and we compute their Gelfand--Kirillov dimensions (Theorems~\ref{D-deformation} and \ref{GK-quantum-group}). We conclude by proving that, if $U$ is an affine prime PI-algebra and a Galois ring, then $\Frac \, \Gamma$ is an algebraic extension of $\Frac \, Z(U)$ (Corollary~\ref{coincide-1}). Moreover, if $U$ is Noetherian, all five dimensions coincide and are equal to $\operatorname{tr.deg}{\k} Z(U)$ (Corollary~\ref{prime-PI-Galois}).

Section~6 is devoted to the affine Hecke algebras constructed in \cite{GKV}. We first show that these algebras are Galois rings (Theorem~\ref{Hecke-1}). When the generalized Cartan matrix is of finite type, the corresponding algebras $H_q(T_{\mathcal C},A)$ are shown to be prime and Noetherian (Theorem~\ref{Hecke-2}). We prove that they are PI-algebras and determine their polynomial identities, PI-degree, and PI-exponent (Theorem~\ref{identities-Hecke}). We further show that, over $\mathbb{C}$, they satisfy the Nullstellensatz and are Jacobson semiprimitive (Theorem~\ref{Nulls-affine-Hecke} and Corollary~\ref{prop-semiprime}). In addition, we prove that they are integral over their centers and therefore satisfy the properties GU, LO, INC, and GD with respect to their centers (Theorem~\ref{Hecke-integral}). We also establish that they satisfy the maximal Nullstellensatz with respect to their centers (Theorem~\ref{Nulls-affine-Hecke}). Finally, we prove that all five ring-theoretic dimensions introduced in Section~5 coincide and are equal to $n$, the size of the square matrix $A$ (Theorem~\ref{dimensions-affine-Hecke}). When $A$ is a generalized Cartan matrix of untwisted affine type, these algebras remain prime (Theorem~\ref{Hecke-4}), and their spherical subalgebras satisfy the Gelfand--Kirillov conjecture (Theorem~\ref{GKC-DAHA}). This completes the picture, since the conjecture for the trigonometric and rational degenerations was previously established in \cite{Schwarz} and \cite{EFOS}, respectively.

The seventh and final section is devoted to the ring-theoretic aspects of spherical Coulomb branch algebras (in the terminology of \cite{lots}). As shown in \cite{Webster} (see also Theorem~\ref{Coulomb-is-Galois-ring}), these algebras are Galois rings. Our results show, in precise terms, that they share many ring-theoretic properties with the universal enveloping algebras of finite-dimensional complex Lie algebras. We prove that, over $\C$, spherical Coulomb branch algebras satisfy the Nullstellensatz and are Jacobson rings (Theorem~\ref{Null+Jacobson}). We also show that they are semiprimitive (Corollary~\ref{CB-semiprimitive}) and are not PI-algebras (Theorem~\ref{Coulomb-not-PI}). Furthermore, we prove that the Gelfand--Tsetlin subalgebra (in the terminology of \cite{lots}) is a Harish--Chandra subalgebra and is maximal among commutative subalgebras with respect to inclusion (Theorem~\ref{Coulomb-maximal-commutative}). Next, we establish the Gelfand--Kirillov conjecture for spherical Coulomb branch algebras (Theorem~\ref{Coulomb-GKC}). An important consequence is that the universal enveloping algebras of the simple Lie algebras of types $B$, $D$, $F$, and $E$ cannot be realized as spherical Coulomb branch algebras, since the Gelfand--Kirillov conjecture fails for them by a theorem of Premet \cite{Premet}.
 We also establish the Gelfand--Kirillov conjecture for the Poisson algebras $\C[\M_C(A,N)]$ in the quasiclassical limit (Theorem~\ref{Coulomb-Poisson}). Finally, we compute the Gelfand--Kirillov and Krull dimensions of spherical Coulomb branch algebras and prove that they are $\LD$-stable (Theorems~\ref{CB-LD-GK} and \ref{Krull-Coulomb}). As a consequence, we obtain the Krull dimension of finite $W$-algebras of type $A$ (Corollary~\ref{Krull-W}). This computation may be viewed as a substantial generalization of Levasseur's theorem on the Krull dimension of universal enveloping algebras of semisimple Lie algebras \cite{Levasseur?}.

\section{Basic definitions and properties}

\subsection{Galois rings}
We work in the framework of \cite{Hartwig}. Namely, let $\Lambda$ be an integrally closed domain, let $G$ be a finite subgroup of $\Autk_\k \, \Lambda$, and let $\mathcal{M}$ be a submonoid of $\Autk_\k \, \Lambda$ satisfying
\begin{enumerate}
    \item $\mathcal{M}\mathcal{M}^{-1}\cap G=\{e\}$;
    \item $G$ acts on $\mathcal{M}$ by conjugation, that is,
    \[
    g\cdot\mu=g\mu g^{-1}, \qquad g\in G,\ \mu\in\mathcal M;
    \]
    \item $\Lambda$ is a Noetherian $\Lambda^G$-module.
\end{enumerate}

The last condition is automatically satisfied when $\Lambda$ is a finitely generated algebra, by Noether's theorem \cite{NoetherX}.

Let $L=\Frac \Lambda$ and consider the skew monoid ring
\[
\mathcal{L}=L*\mathcal{M}.
\]
Set
\[
\Gamma=\Lambda^G,\qquad K=\Frac\Gamma.
\]
Then $K=L^G$. We also define
$
\mathcal K=\mathcal L^G,
$
where the action of $G$ on $L*\mathcal M$ is given by
\[
(a\mu)^g=g(a)\,(g\cdot\mu),
\qquad
a\in L,\ \mu\in\mathcal M.
\]

\begin{proposition}
\begin{enumerate}[{\rm (i)}]
    \item $\Lambda$ is integral over $\Gamma$.
    \item $\Gamma$ is integrally closed.
    \item $\Lambda$ is the integral closure of $\Gamma$ in $L$.
    \item $\Lambda$ is a finitely generated $\Gamma$-module and a Noetherian ring.
    \item $\Gamma$ is a Noetherian ring.
\end{enumerate}
\end{proposition}

\begin{proof}
\cite[Lemma~2.1]{Hartwig}.
\end{proof}

\begin{definition}\cite{DFO}
Let $U$ be a ring and $C$ a commutative subring. We say that $C$ is a Harish--Chandra subring if, for every $u\in U$, the bimodule $CuC$ is finitely generated as both a left and a right $C$-module.
\end{definition}

Further studies of Harish--Chandra subrings and their generalizations can be found in \cite{Fillmore,Schwarz2}.

\begin{definition}\cite{FO}
\begin{itemize}
\item Let $U$ be a finitely generated $\Gamma$-subring of $\mathcal K$. Then $U$ is called a \emph{Galois $\Gamma$-ring} if
\[
KU=UK=\mathcal K.
\]

\item Let $U$ be a Galois $\Gamma$-ring in $\mathcal K$. If, for every finite-dimensional left (respectively, right) $K$-subspace $W$ of $\mathcal K$, the intersection $W\cap U$ is a finitely generated left (respectively, right) $\Gamma$-module, then $U$ is called a \emph{Galois $\Gamma$-order}.
\end{itemize}
\end{definition}

\begin{theorem}\label{HC-subalgebra}
$\Gamma$ is a Harish--Chandra subring of every Galois $\Gamma$-ring containing it. Moreover, if $U$ is a Galois $\Gamma$-order, then $\Gamma$ is maximal among commutative subalgebras of $U$.
\end{theorem}

\begin{proof}
\cite[Lemma~2.4, Proposition~2.4, Proposition~2.14]{Hartwig}.
\end{proof}

\textbf{Convention.}
Throughout the paper we assume that $\mathcal M$ has finitely many $G$-orbits. In particular, $\mathcal M$ is finitely generated as a monoid. For $\mu\in\mathcal M$, we denote its orbit under the action of $G$ by $\mathcal O_\mu$.

\subsection{Generalized Weyl algebras}




We recall the definition of a \emph{generalized Weyl algebra} (abbreviated GWA), introduced by V. Bavula \cite{Bavula}.

\begin{definition}
Let $D$ be a ring, let $\sigma=(\sigma_1,\ldots,\sigma_n)$ be an $n$-tuple of commuting automorphisms of $D$, and let $a=(a_1,\ldots,a_n)$ be an $n$-tuple of nonzero central elements of $D$ satisfying
\[
\sigma_i(a_j)=a_j,\qquad i\neq j.
\]
The \emph{generalized Weyl algebra} $D(a,\sigma)$ of rank $n$ is the algebra generated over $D$ by elements
\[
X_i^+,\,X_i^-,\qquad i=1,\ldots,n,
\]
subject to the relations
\begin{subequations}\label{eq:GWA-relations}
\begin{gather}
X_i^+d=\sigma_i(d)X_i^+,\qquad
X_i^-d=\sigma_i^{-1}(d)X_i^-,
\quad \forall\, d\in D,\\
[X_i^+,X_j^+]=[X_i^-,X_j^-]=[X_i^+,X_j^-]=0,
\qquad i\neq j,\\
X_i^-X_i^+=a_i,\qquad
X_i^+X_i^-=\sigma_i(a_i).
\end{gather}
\end{subequations}
\end{definition}

The algebra $D(a,\sigma)$ is free as both a left and a right $D$-module. Moreover, if $D$ is a domain, then $D(a,\sigma)$ is also a domain. If $D$ is left or right Noetherian, then so is $D(a,\sigma)$ \cite{Bavula}. An excellent survey of generalized Weyl algebras and related constructions is \cite{Gaddis}.

\begin{example}\label{Witten-Woronowicz}
The Witten--Woronowicz deformation is an example of a generalized Weyl algebra \cite{BO}. It is given by $D(a,\sigma)$, where
\[
D=\mathbb{C}[H,Z],\qquad
a=Z+\alpha H+\beta,
\]
and
\begin{gather*}
\sigma(H)=s^4H,\qquad
\sigma(Z)=s^2Z,\\
\alpha=-\frac{1}{s(1-s^2)},\qquad
\beta=\frac{s}{1-s^4},
\end{gather*}
where $s\in\mathbb C^\times$ satisfies $s^4\neq1$.
\end{example}

Let $D(a,\sigma)$ be a generalized Weyl algebra of rank $n$. Since $D(a,\sigma)$ is generated as a $D$-algebra by elements $\{x_i\}_{i\in I}$ satisfying
\[
Dx_i=x_iD,\qquad i\in I,
\]
it follows immediately that $D$ is a Harish--Chandra subring of $D(a,\sigma)$.

To fit into our framework, we now assume that $D$ is a commutative integrally closed Noetherian domain. Let $\mathcal M$ denote the subgroup of $\Autk \, D$ generated by $\sigma_1,\ldots,\sigma_n$. If the natural homomorphism
\[
\mathbb Z^n\longrightarrow\mathcal M
\]
is an isomorphism, then, by a straightforward generalization of \cite[Theorem~14]{FS3} with the same proof, we obtain the following result.

\begin{proposition}\label{GWA}
Let $F=\Frac D$ and let $\mathcal K=F*\mathcal M$. Then $D(a,\sigma)$ is a Galois $D$-ring in $\mathcal K$.
\end{proposition}


\subsection{The center of Galois rings}

We now study the center of a Galois ring. Our goal is to obtain an analogue of \cite[Theorem~4.1(4)]{FO}.

For $\mu\in\mathcal M$, let $G_\mu$ denote its stabilizer under the action of $G$.

\begin{definition}
Let $\mu\in\mathcal M$ and let $a\in L^{G_\mu}$. We define
\[
[a\mu]=\sum_{g\in G/G_\mu} g(a)\,g\cdot\mu.
\]
Then $[a\mu]$ is an element of $\mathcal K$.
\end{definition}

It is clear that this definition is independent of the choice of representatives of the cosets in $G/G_\mu$. We have the following analogue of \cite[Lemma~2.1]{FO}.

\begin{proposition}\label{FO}
\begin{enumerate}[{\rm (i)}]
\item The element $[a\mu]$ is $G$-invariant. Consequently,
$
[a\mu]\in\mathcal K.
$

\item Let $\mu\in\mathcal M$ and define
\[
\mathcal K_\mu=\{[a\mu]\mid a\in L^{G_\mu}\}.
\]
Then $\mathcal K_\mu$ is an $L^{G_\mu}$-bimodule, and hence a $K=L^G$-bimodule. Moreover, for every $\gamma\in L^{G_\mu}$,
\[
\gamma[a\mu]=[\gamma a\,\mu],
\qquad
[a\mu]\gamma=[a\,\mu(\gamma)\,\mu].
\]

\item If
\[
\mathcal M=\mathcal O_{\mu_1}\sqcup\cdots\sqcup\mathcal O_{\mu_r}
\]
is the decomposition of $\mathcal M$ into $G$-orbits, then
\[
\mathcal K=\bigoplus_{i=1}^r\mathcal K_{\mu_i}
\]
as a $K$-bimodule.
\end{enumerate}
\end{proposition}

\begin{proof}
The proof is identical to that of \cite[Lemma~2.1]{FO}.
\end{proof}

Our next objective is to extend \cite[Theorem~4.1]{FO} to the framework of \cite{Hartwig}.

\begin{theorem}\label{center}
Let $U$ be a $\Gamma$-ring contained in $(L*\mathcal M)^G$.

\begin{enumerate}[{\rm (i)}]

\item $U\cap K$ is a maximal commutative subalgebra of $U$.

\item
\[
Z(U)=(U\cap K)^{\mathcal M}.
\]

\item
If $U$ is a Galois $\Gamma$-order in $\mathcal K$, then
\[
U\cap K=\Gamma,
\]
and hence
\[
Z(U)=\Gamma^{\mathcal M}.
\]
\end{enumerate}
\end{theorem}

\begin{proof}
Statements (i) and (ii) are proved \emph{mutatis mutandis} by the same arguments as the corresponding statements in \cite[Theorem~4.1]{FO}. Statement (iii) follows from \cite[Proposition~2.14]{Hartwig}.
\end{proof}

As a simple illustration of the method, we obtain the following consequence.

\begin{proposition}
The center of the Witten--Woronowicz deformation is $\mathbb C$ whenever $s$ is not a root of unity.
\end{proposition}

\section{Localization and applications of Goldie's Theorem}

One of the main reasons for the success of the theory of Galois rings in computing skew fields of fractions is the following observation. If $U$ is a Galois ring in $\mathcal K$ and both $U$ and $\mathcal K$ are Ore domains, then they have the same skew field of fractions. This follows from the more precise statement \cite[Proposition~4.2]{FO}, which we now adapt to the framework of arbitrary rings, rather than algebras over a field.

\begin{theorem}\label{localization2}
Let $S=\Gamma\setminus\{0\}$. Then every element of $S$ is regular, $S$ is both a left and a right Ore set, and the left and right localizations of the Galois ring $U$ with respect to $S$ are both isomorphic to $\mathcal K$.
\end{theorem}

\begin{proof}
For the first claim, let $\gamma\in S$ and let $u\in U$ be nonzero. We must show that $\gamma u\neq0$ and $u\gamma\neq0$. Since $U\hookrightarrow\mathcal K$, we may write
\[
u=\sum_{\mu}\ell_\mu\mu,
\]
where only finitely many coefficients $\ell_\mu\in L$ are nonzero. Then
\[
\gamma u=\sum_{\mu}\gamma\ell_\mu\,\mu,
\]
and each coefficient $\gamma\ell_\mu$ is nonzero because $\gamma\neq0$ in the domain $L$. Similarly,
\[
u\gamma=\sum_{\mu}\ell_\mu\,\mu(\gamma)\,\mu,
\]
and $\mu(\gamma)\neq0$ since $\Gamma\subset K\subset L$.

We now prove the left Ore condition and show that $S^{-1}U\simeq\mathcal K$. The right-handed argument is symmetric. By \cite[Lemma~2.7(iv)]{Hartwig}, we have
\[
KU=UK=\mathcal K,
\qquad
K=\operatorname{Frac}\Gamma.
\]
Let $s\in S$ and $u\in U$. Since $s^{-1}u\in\mathcal K=UK$, there exist elements $u_i\in U$, $\gamma_i\in\Gamma$, and nonzero elements $s_i\in S$ such that
\[
s^{-1}u=\sum_{i=1}^n u_i\gamma_i s_i^{-1}.
\]
Set
\[
r=s_1s_2\cdots s_n\in S.
\]
After clearing denominators, we obtain
\[
s^{-1}u=u'r^{-1}
\]
for some $u'\in U$. This proves the left Ore condition. The right Ore condition is proved in the same way.

Finally, consider the left annihilator
\[
\operatorname{ass}(S)=\{u\in U\mid (\exists\,s\in S)\ su=0\}.
\]
Since $U\hookrightarrow\mathcal K$ and every $s\in S$ becomes invertible in $\mathcal K$, the equality $su=0$ implies $u=0$. Hence $\operatorname{ass}(S)=(0)$. By \cite[§2.1.3]{McConnell}, the Ore localization of $U$ with respect to $S$ exists and is isomorphic to $\mathcal K$. The same argument applies to the right localization.
\end{proof}

Throughout the paper, a noncommutative Noetherian ring means a ring that is both left and right Noetherian. The following corollary is the noncommutative analogue of a classical result from algebraic geometry.

\begin{corollary}\label{collary-spec}
Assume that $U$ is Noetherian. Then there is a bijection between
\[
\{P\in\operatorname{Spec}U\mid P\cap S=\varnothing\}
\]
and $\operatorname{Spec}\mathcal K$. The correspondence sends a prime ideal $P$ of $U$ to the prime ideal $P\mathcal K$ of $\mathcal K$, and a prime ideal $Q$ of $\mathcal K$ to $Q\cap U$. Moreover, $\mathcal K$ is Noetherian.
\end{corollary}

\begin{proof}
The bijection follows from \cite[Proposition~2.1.16(vii)]{McConnell}, and the Noetherianity of $\mathcal K$ follows from \cite[Corollary~10.16(a)]{GW}.
\end{proof}

In the following definition we introduce two algebras that frequently arise as the ambient skew monoid ring $\mathcal L$ in the definition of a Galois ring. In fact, almost all known examples of Galois rings have this form; see \cite{Hartwig}.

\begin{definition}\label{main-chars}
Assume that $\operatorname{char}\k=0$.

\begin{enumerate}
\item For $m\le n$, let
\[
S_m^n(\k)=\k(x_1,\ldots,x_n)*\mathbb Z^m,
\]
where the canonical basis $\varepsilon_1,\ldots,\varepsilon_m$ of $\mathbb Z^m$ acts by
\[
\varepsilon_i(x_j)=
\begin{cases}
x_j-\delta_{ij}, & j\le m,\\
x_j, & j>m.
\end{cases}
\]

\item For $m\le n$, let
\[
Q_m^n(\k)=\k(x_1,\ldots,x_n)*\mathbb Z^m,
\]
where
\[
\varepsilon_i(x_j)=
\begin{cases}
q^{\delta_{ij}}x_j, & j\le m,\\
x_j, & j>m,
\end{cases}
\]
and $q\in\k^\times$ is not a root of unity.
\end{enumerate}
\end{definition}

\begin{proposition}\label{main-lemma}
The algebras $S_m^n(\k)$ and $Q_m^n(\k)$ are Noetherian simple rings.
\end{proposition}

\begin{proof}
The algebra $S_m^n(\k)$ is a localization of the Weyl algebra, while $Q_m^n(\k)$ is a localization of the quantum torus. The proofs from \cite{FO,FH} carry over verbatim using Theorem~\ref{localization2}. Since the Weyl algebra and the quantum torus are Noetherian simple rings \cite{GW,Brown}, the same holds for $S_m^n(\k)$ and $Q_m^n(\k)$.
\end{proof}

Assume that $\operatorname{char}\k=0$, and let $V$ be an $n$-dimensional $\k$-vector space. Recall that a \emph{$\k$-lattice} in $V$ is a free abelian group $\Sigma$ of rank $n$ such that
\[
V=\k\otimes_{\mathbb Z}\Sigma.
\]
The lattice $\Sigma$ acts on $V$ by translations,
\[
v\longmapsto v+\lambda,
\qquad
v\in V,\ \lambda\in\Sigma.
\]
Consequently, $\Sigma$ acts on the rational function field
\[
\k(V)=\operatorname{Frac}S(V^*)
\]
via
\[
(\lambda\cdot f)(v)=f(v-\lambda),
\qquad
v\in V,\ \lambda\in\Sigma,\ f\in\k(V).
\]
The following lemma shows that shift operator algebras occur naturally in this context.

\begin{lemma}\label{fundamental-lattices}
There is an isomorphism
\[
\k(V)*\Sigma
\simeq
\k(x_1,\ldots,x_n)*\mathbb Z^n,
\]
where
\[
\varepsilon_i(x_j)=x_j-\delta_{ij}
\qquad
\text{for all } i,j.
\]

Moreover, if $W<\operatorname{GL}(\Sigma)$ is a finite subgroup acting on $\Sigma$ by integral linear transformations, then
\[
(\k(V)*\Sigma)^W
\simeq
(\k(x_1,\ldots,x_n)*\mathbb Z^n)^W,
\]
where $W$ acts linearly on the variables $x_i$ and normalizes $\mathbb Z^n$ by conjugation.
\end{lemma}

\begin{proof}
Choose a $\mathbb Z$-basis $v_1,\ldots,v_n$ of $\Sigma$. Since
\[
V=\k\otimes_{\mathbb Z}\Sigma,
\]
the vectors $v_1,\ldots,v_n$ also form a basis of $V$. Let $x_1,\ldots,x_n\in V^*$ be the dual basis. In these coordinates, the action of $\Sigma$ is exactly the standard translation action of $\mathbb Z^n$, and the statement follows immediately.
\end{proof}

\begin{proposition}
If $U$ is simple, then
\[
Z(U)=Z(\mathcal K).
\]
\end{proposition}

\begin{proof}
This follows from Theorem~\ref{localization2} and \cite[Proposition~2.1.16(viii)]{McConnell}.
\end{proof}

\begin{corollary}\label{to-use-PI}
\[
Z\bigl(S_m^n(\k)\bigr)
=
Z\bigl(Q_m^n(\k)\bigr)
=
\k(z_1,\ldots,z_{n-m}),
\]
where $z_i=x_{m+i}$.
\end{corollary}

\begin{proof}
It is well known that the centers of the Weyl algebra and the quantum torus consist only of scalars; see \cite{Gelfand,Brown}. Since these algebras are simple, the claim follows from the previous proposition.
\end{proof}

Because of their central role in the theory of Galois rings, the algebras $S_m^n(\k)$ and $Q_m^n(\k)$ will appear repeatedly throughout the remainder of the paper.



\subsection{A sufficient condition for \texorpdfstring{$U$}{U} to be Ore}

The objective of this subsection is to prove that if $\mathcal{M}$ is torsion-free and locally polycyclic-by-finite, then every Galois ring $U$ is an Ore domain (Theorem~\ref{U-is-Ore}).

\begin{definition}
Let $R$ be a ring, and suppose that the set $S$ of regular elements is a right denominator set. We denote by $\Q(R)$ the right localization of $R$ at $S$ and call it the \emph{classical right quotient ring} (or \emph{classical right ring of fractions}) of $R$.
\end{definition}

When $R$ is an Ore domain, we write $\Frac ], R$ instead of $\Q(R)$ and refer to it as the \emph{skew field of fractions} of $R$.

We recall the notion of a right Goldie ring.

\begin{definition}
A ring $R$ is called a \emph{right Goldie ring} if it has finite right uniform dimension and satisfies the ascending chain condition on right annihilators.
\end{definition}

In particular, every right Noetherian ring is right Goldie. For the reader's convenience, we recall Goldie's theorem in the form needed below.

\begin{theorem}\label{Goldie}
\begin{enumerate}[{\rm (i)}]
\item A ring $R$ has a classical right quotient ring that is semisimple Artinian if and only if $R$ is a semiprime right Goldie ring.

\item A ring $R$ has a classical right quotient ring that is simple Artinian if and only if $R$ is a prime right Goldie ring.
\end{enumerate}
\end{theorem}

\begin{proof}
See \cite[Theorem~2.3.6]{McConnell}.
\end{proof}

\begin{theorem}\label{main}
Let $U$ be a Galois $\Gamma$-ring in $\mathcal K$. Then $\Q(U)$ exists if and only if $\Q(\mathcal K)$ exists. In this case,
\[
\Q(U)\cong\Q(\mathcal K).
\]
\end{theorem}

\begin{proof}
Suppose first that $\Q(U)$ exists. Since localization is transitive, Theorem~\ref{localization2} implies that $\Q(\mathcal K)$ also exists, and the two quotient rings are canonically isomorphic.

Conversely, suppose that $\Q(\mathcal K)$ exists. Let
\[
x=ab^{-1}\in\Q(\mathcal K),
\]
where $b$ is regular in $\mathcal K$. By Theorem~\ref{localization2},
\[
S=\Gamma\setminus\{0\}
\]
is a left and right Ore set consisting of regular elements, and every element of $\mathcal K$ can be written as a fraction over $S$. Hence
\[
a=us^{-1},
\qquad
b=vr^{-1},
\]
for some $u,v\in U$ and $r,s\in S$. Since $b$ is regular in $\mathcal K$, the element $v$ is regular in $U$.

Because $S$ is a left Ore set in $U$, there exist $w\in U$ and $t\in S$ such that
\[
tv=sw.
\]
Therefore,
\[
wt^{-1}=s^{-1}v
\]
in $\Q(\mathcal K)$, and consequently
\[
ab^{-1}
=
uw\,s^{-1}t^{-1}
=
(uw)(ts)^{-1},
\]
since $s^{-1}$ and $t^{-1}$ belong to the commutative field $K=\Frac \, \Gamma$. Thus every element of $\Q(\mathcal K)$ belongs to the classical quotient ring of $U$, and the conclusion follows from \cite[2.1.3]{McConnell}.
\end{proof}

\begin{theorem}\label{localization4}
The ring $U$ is prime right Goldie if and only if $\mathcal K$ is prime right Goldie. Likewise, $U$ is semiprime right Goldie if and only if $\mathcal K$ is semiprime right Goldie.
\end{theorem}

\begin{proof}
By Theorem~\ref{main}, the classical right quotient ring of $U$ exists if and only if that of $\mathcal K$ exists, and in this case they are canonically isomorphic. By Goldie's theorem, $U$ is prime right Goldie (respectively, semiprime right Goldie) if and only if its classical quotient ring is simple Artinian (respectively, semisimple Artinian). The same characterization holds for $\mathcal K$, and the result follows.
\end{proof}

If the Goldie assumption is omitted, we can still transfer primeness and semiprimeness from $U$ to $\mathcal{K}$.

\begin{proposition}\label{prop-semiprime}
If $U$ is a prime (respectively, semiprime) ring, then so is $\mathcal{K}$.
\end{proposition}

\begin{proof}
See \cite[Proposition~3.1.15]{Rowen}.
\end{proof}

We also have an analogue of Theorem~\ref{main}. Its proof is identical, using the well-known characterization that a domain has a skew field of fractions if and only if it is an Ore domain.

\begin{theorem}\label{Ore}
$U$ is an Ore domain if and only if $\mathcal{K}$ is an Ore domain.
\end{theorem}

The theorem above allows us to reduce the problem to $\mathcal{K}$. We now go one step further and obtain a general condition on the monoid $\mathcal{M}$ that guarantees that $U$ is an Ore domain. This explains why the Ore property holds in virtually all known examples of Galois rings (cf.~\cite[Introduction]{Schwarz}).

If $D$ is an integral domain and $G<\operatorname{Aut}(D)$ is a finite group of ring automorphisms, then it is an elementary exercise to prove that
\[
\Frac(D^G)=\Frac(D)^G.
\]
The corresponding statement in the noncommutative setting is considerably deeper.

\begin{theorem}\label{Faith}
Let $A$ be an Ore domain, and let $G$ be a finite group acting faithfully on $A$ by ring automorphisms. Then $A^G$ is an Ore domain and
\[
\Frac(A^G)=\Frac(A)^G.
\]
\end{theorem}

\begin{proof}
See \cite{Faith}.
\end{proof}

The next result shows that, in order to prove that $U$ is an Ore domain, it suffices to verify that the skew monoid ring $\mathcal{L}=L*\mathcal{M}$ is an Ore domain. This reduces the problem entirely to the structure of the monoid $\mathcal{M}$.

\begin{proposition}\label{localization3}
If $\mathcal{L}$ is a left (respectively, right) Ore domain, then so is $U$.
\end{proposition}

\begin{proof}
By Theorem~\ref{Faith}, the fixed ring
\[
\mathcal{L}^G=\mathcal{K}
\]
is a left (respectively, right) Ore domain. The conclusion now follows from Theorem~\ref{Ore}.
\end{proof}

We now obtain a criterion, depending only on the monoid $\mathcal{M}$, that guarantees that every Galois ring is an Ore domain.

\begin{definition}
A group $G$ is called \emph{polycyclic-by-finite} if it admits a normal series
\[
\{1\}=G_0\triangleleft G_1\triangleleft\cdots\triangleleft G_n=G,
\]
whose successive quotients are either finite or infinite cyclic. A group is called \emph{locally polycyclic-by-finite} if every finitely generated subgroup is polycyclic-by-finite.
\end{definition}

It is well known (see \cite[1.5.12]{McConnell}) that the crossed product of a Noetherian ring with a polycyclic-by-finite group is again Noetherian. An important open problem in the theory of group rings asks whether the converse holds, namely, whether the Noetherianity of the integral group ring $\mathbb{Z}G$ implies that $G$ is polycyclic-by-finite (see \cite{Passman2}, \cite[§8.2]{Rowen2}).

The following theorem, due to Passman \cite[Theorem~37.10]{Passman2}, will play a crucial role.

\begin{theorem}\label{Passman}
Let $D$ be an Ore domain, and let $G$ be a torsion-free\footnote{Important examples of torsion-free groups are totally ordered groups.} locally polycyclic-by-finite group acting faithfully on $D$ by ring automorphisms. Then the skew group ring $D*G$ is an Ore domain. In particular, $D*\mathbb{Z}^n$ is an Ore domain.
\end{theorem}

Combining the preceding results, we obtain the following theorem.

\begin{theorem}\label{U-is-Ore}
Let $U$ be a Galois $\Gamma$-ring. If $\mathcal{M}$ is a torsion-free locally polycyclic-by-finite group, then $U$ is an Ore domain. In particular, if $\mathcal{M}\cong\mathbb{Z}^n$, then $U$ is an Ore domain.
\end{theorem}

\begin{proof}
This follows immediately from Proposition~\ref{localization3} and Theorems~\ref{Passman} and \ref{Faith}.
\end{proof}

\subsection{Prime Goldie Galois rings} 

In this subsection we show that, under natural assumptions, Galois rings are prime right Goldie rings.

\begin{definition}
Let $R$ be a prime right Goldie ring. An automorphism $\phi$ of $R$ is called \emph{X-outer}\footnote{See \cite[Example~3.7]{Montgomery}.} if there is no non-zero element $x$ of $\Q(R)$ such that $x \phi(r)=r x, \forall r \in R$. A group $G$ of automorphisms of $R$ is called \emph{X-outer} if every nontrivial element of $G$ is X-outer.
\end{definition}

In case $R$ is an Ore domain, an automorphism $\phi$ then is X-outer if its restriction to $\Frac \, R$ is an outer automorphism. If $R$ is a field every non-trivial automorphism is X-outer.

\begin{theorem}\label{proposition-Fisher-Montgomery}
Let $R$ be a prime Goldie ring, and let $G$ be a group acting faithfully on $R$ by X-outer automorphisms. Then the skew group ring $R*G$ is prime Goldie. If $R$ is simple, then $R*G$ is simple. Moreover, if $G$ is polycyclic-by-finite, then $R*G$ is Noetherian.
\end{theorem}

\begin{proof}
The first two assertions follow from \cite[Corollary~3]{FisherM}. As noted above, if $G$ is polycyclic-by-finite, then $R*G$ is Noetherian.
\end{proof}

Before proceeding, we recall the following fundamental results from noncommutative invariant theory.

Observe that, in the setting of Galois rings, if $\mathcal M$ is nontrivial, then its action on the field $L$ is automatically X-outer.

\begin{theorem}\label{invariant-theory}
Let $R$ be a simple algebra, and let $G$ be a finite group of outer automorphisms of $R$ whose order is invertible in the base field. Then
\begin{enumerate}[{\rm (i)}]
\item $R^G$ is simple.
\item $R$ is a finitely generated projective $R^G$-module.
\item $R^G$ and $R*G$ are Morita equivalent.
\end{enumerate}
\end{theorem}

\begin{proof}
See \cite[Lemma~2.1, Theorems~2.4 and~2.5, Corollary~2.6]{Montgomery}.
\end{proof}

\begin{theorem}\label{U-is-prime-Goldie}
Let $U$ be a Galois $\Gamma$-ring in $\mathcal K=\mathcal L^G$, where $\mathcal M$ is a nontrivial group.

\begin{enumerate}[{\rm (i)}]
\item The skew group ring $\mathcal L=L*\mathcal M$ is simple.

\item If $G$ acts on $\mathcal L$ by outer automorphisms and $|G|^{-1}\in L$, then $\mathcal K$ is simple.

\item If $G$ acts on $\mathcal L$ by X-outer automorphisms and $\mathcal M$ is polycyclic-by-finite, then $\mathcal K$ is a prime Noetherian ring. Consequently, $U$ is a prime right Goldie ring.
\end{enumerate}
\end{theorem}

\begin{proof}
Since $L$ is a field, it is a simple ring. Therefore, Theorem~\ref{proposition-Fisher-Montgomery} implies that
\[
\mathcal L=L*\mathcal M
\]
is simple.

If $G$ acts on $\mathcal L$ by outer automorphisms and $|G|^{-1}\in L$, then Theorem~\ref{invariant-theory} implies that
$
\mathcal K=\mathcal L^G
$
is simple.

Finally, if $G$ acts by X-outer automorphisms and $\mathcal M$ is polycyclic-by-finite, then $\mathcal L$ is Noetherian by Theorem~\ref{proposition-Fisher-Montgomery}. Hence \cite[Corollary~1.12 and Theorem~3.17]{Montgomery} imply that $\mathcal K$ is prime and Noetherian. In particular, $\mathcal K$ is prime right Goldie, and therefore Theorem~\ref{localization4} yields that $U$ is prime right Goldie.
\end{proof}

With this we have a rather remarkable consequence.

\begin{corollary}
    Every prime ideal of $U(\mathfrak{gl}_n)$ intersects the Gelfand-Tsetlin subalgebra $\mathcal{Z}$ non-trivially.
\end{corollary}
\begin{proof}
    This is an immediate consequence of the Galois ring realization of $U(\mathfrak{gl}_n)$ in \cite{FO}, Theorem \ref{U-is-prime-Goldie}(ii), and Corollary \ref{collary-spec}.
\end{proof}

\section{PI Galois rings}

In this section we study PI Galois rings. By combining the powerful techniques available in the theories of PI-algebras and Galois rings, we obtain several structural results. We begin by recalling the notion of a PI-algebra.

\begin{definition}
Let $X={x_1,x_2,\ldots}$ be a countable set, let $\Phi$ be a commutative unital ring, and let $\Phi\langle X\rangle$ denote the free associative unital $\Phi$-algebra generated by $X$. Let
$
f(x_1,\ldots,x_n)\in\Phi\langle X\rangle
$
be a noncommutative polynomial involving only the variables $x_1,\ldots,x_n$, and assume that at least one monomial of maximal degree has coefficient $1$.\footnote{This assumption excludes trivial identities such as the polynomial $px$, which vanishes identically in every ring of characteristic $p$.} A $\Phi$-algebra $A$ is called a \emph{PI-algebra} if there exists such a polynomial $f$ satisfying
$
f(a_1,\ldots,a_n)=0
$
for all $a_1,\ldots,a_n\in A$. In this case, $f$ is called a \emph{polynomial identity} of $A$. If $\Phi=\mathbb{Z}$, we simply say that $A$ is a \emph{PI-ring}.
\end{definition}

\begin{definition}
Let $A$ be a PI-algebra over a commutative ring $\Phi$. The \emph{$T$-ideal} of $A$, denoted by $T(A)$, is the ideal of $\Phi\langle X\rangle$ consisting of all polynomial identities of $A$.
\end{definition}

If $\Phi$ is a field $\k$ of characteristic zero, the notions of PI-ring and PI-algebra coincide by a theorem of Amitsur (see \cite[Part~B, Chapter~7]{Drensky}). The converse implication is immediate: every PI-ring that is a $\k$-algebra is automatically a PI-algebra. For further background on PI-algebras we refer the reader to \cite{Drensky}, \cite{GZ}, \cite[Chapter~13]{McConnell}, and \cite[Chapter~6]{Rowen2}.

We shall repeatedly use the following elementary fact.

\begin{lemma}\label{baby}
Let $A$ be a PI-algebra, $B$ a homomorphic image of $A$, and $C$ a subalgebra of $A$. Then both $B$ and $C$ are PI-algebras, and
$
T(A)\subseteq T(B),
\,\,
T(A)\subseteq T(C).
$
\end{lemma}

\begin{theorem}\label{stable}
Let $A$ be a PI-algebra over an infinite integral domain $\Phi$, and suppose that $A$ is torsion-free as a $\Phi$-module. Then
$
T(A)=T(A\otimes_\Phi H)
$
for every torsion-free commutative $\Phi$-algebra $H$.
\end{theorem}

\begin{proof}
See \cite[Theorem~6.1.44]{Rowen2}.
\end{proof}

We shall also use Posner's Theorem.

\begin{theorem}\label{Posner}
Let $R$ be a prime PI-ring with center $Z$, and let $Q$ be the localization of $R$ at the nonzero elements of $Z$. Let $F=\Frac \, Z$. Then $Q$ is a finite-dimensional central simple $F$-algebra, it is the classical quotient ring of $R$, and $Q$ and $R$ satisfy the same polynomial identities in $\mathbb{Z}\langle X\rangle$. If, moreover, $R$ is a $\k$-algebra over a field $\k$, then
$
T(R)=T(Q)=T(M_n(\k)),
$
where
$
n=\sqrt{[Q]}.
$
\end{theorem}

\begin{proof}
See \cite[Theorem~B.6.5]{Drensky} and \cite[Theorem~1.11.13]{GZ}.
\end{proof}

\begin{definition}
In the setting of Posner's Theorem, the integer
$
n=\sqrt{[Q]}
$
is called the \emph{PI-degree} of $R$ and is denoted by
$
\operatorname{PI-deg}(R).
$
\end{definition}

\begin{theorem}\label{PI-0}
Let $U$ be a Galois algebra over an infinite field $\k$. Then $U$ is a PI-algebra if and only if $\K$ is a PI-algebra. Moreover,
\[
T(U)=T(\K).
\]
\end{theorem}

\begin{proof}
Let $A:=U\otimes_{\k}K$. By Theorem \ref{stable},
\[
T(U)=T(A).
\]
By the definition of a Galois ring, $\K$ is a homomorphic image of $A$. Hence,
\[
T(U)=T(A)\subseteq T(\K).
\]
On the other hand, since $U\subseteq\K$, we have
\[
T(\K)\subseteq T(U).
\]
Therefore,
\[
T(U)=T(\K),
\]
as claimed.
\end{proof}

\begin{corollary}\label{not-PI}
Suppose that $U$ is a PI-algebra without nilpotent elements. Then $\L$ is also a PI-algebra. The same conclusion holds if $U$ is a semiprime PI-algebra over a commutative ring $\Phi$ of characteristic $0$.
\end{corollary}

\begin{proof}
The first statement follows from Theorem \ref{PI-0} together with \cite[Theorem 6.8]{Montgomery}. The second follows from Theorem \ref{PI-0}, Proposition \ref{prop-semiprime}, and \cite[Theorem 6.5]{Montgomery}.
\end{proof}

\begin{theorem}\label{usually-not-PI}
Suppose that $U$ is a Galois ring in
\[
\K=\L^G=(L*\M)^G,
\]
where $\L$ is either $Q_m^n(\k)$ or $S_m^n(\k)$. Then $U$ is not a PI-algebra.
\end{theorem}

\begin{proof}
By Theorem \ref{U-is-Ore}, $U$ is an Ore domain. Hence, by the previous corollary, if $U$ were a PI-algebra, then $\L$ would also be a PI-algebra. However, $\L$ is either $Q_m^n(\k)$ or $S_m^n(\k)$, both of which are simple and infinite-dimensional over their centers (Proposition \ref{main-lemma} and Corollary \ref{to-use-PI}). Therefore, by Kaplansky's theorem (\cite[Theorem 1.11.7]{GZ}), neither of these algebras is PI. This contradiction shows that $U$ cannot be a PI-algebra.
\end{proof}

With this result, we obtain the following theorem.

\begin{theorem}\label{list}
The following algebras are not PI-algebras:
\begin{enumerate}[{\rm (i)}]
\item
The alternating analogue $\mathfrak{A}(\mathfrak{gl}_n)$ of $U(\mathfrak{gl}_n)$ \cite{Jauch}.
\item
Finite $W$-algebras of type $A$, OGZ algebras, quantum OGZ algebras, and their parabolic subalgebras \cite{Hartwig}. In particular, $U_q(\mathfrak{gl}_n)$.
\item
Spherical subalgebras of rational Cherednik algebras associated with the complex reflection groups $G(m,p,n)$ \cite{LW}\footnote{When the deformation parameter $\mathfrak{c}$ is not zero.}.
\item
Spherical subalgebras of trigonometric Cherednik algebras associated with $S_n$ \cite{KN}.
\item
Deformations of Kleinian singularities of types $A$ and $D$ \cite{Hartwig3}.
\end{enumerate}
\end{theorem}

Item~(i) of the preceding theorem is new. The remaining algebras were already known not to be PI-algebras by other methods. Nevertheless, this theorem illustrates the effectiveness of realizing an algebra as a Galois ring in proving that it is not a PI-algebra. Later, we will show that spherical Coulomb branch algebras are also not PI-algebras (Theorem~\ref{Coulomb-not-PI}).

\begin{theorem}\label{PI}
Let $U$ be a Galois $\Gamma$-ring in the invariant skew group ring
\[
\mathcal{K}=\mathcal{L}^G,\qquad \mathcal{L}=L*\mathcal{M}.
\]
Then $U$ is a prime PI $\k$-algebra if and only if $\mathcal{K}$ is a central simple algebra. In this case,
\[
\Frac(Z(U))=Z(\mathcal{K}),
\]
the group $\mathcal{M}$ is necessarily finite, and $U$ and $\mathcal{K}$ satisfy the same polynomial identities as the matrix algebra
\[
M_n(\Frac(Z(U))),
\]
where
\[
n=\operatorname{PI-deg}(U).
\]
\end{theorem}

\begin{proof}
Suppose first that $\mathcal{K}$ is a central simple algebra. Then $\mathcal{K}$ is, in particular, a prime Goldie ring. By Proposition~\ref{localization4}, it follows that $U$ is also a prime ring. Moreover, since $\mathcal{K}$ is a PI-ring, so is $U$ by Lemma~\ref{baby}.

Conversely, assume that $U$ is a prime PI-ring. Let $Z=Z(U)$ denote its center. By Posner's theorem,
\[
\Q(U)=UZ^{-1}
\]
is a finite-dimensional central simple algebra over $\Frac(Z)$. Since
\[
\Q(U)=UZ^{-1}\subseteq U\Gamma^{-1}=\mathcal{K}\subseteq \Q(U),
\]
we conclude that
\[
\mathcal{K}=\Q(U).
\]
Hence $\mathcal{K}$ is a finite-dimensional central simple algebra over $\Frac(Z(U))$, and therefore
\[
Z(\mathcal{K})=\Frac(Z(U)).
\]
The remaining assertions, including the finiteness of $\mathcal{M}$ and the description of the polynomial identities, follow immediately from Posner's theorem.
\end{proof}

This theorem allows us to compute the PI-exponent of such Galois rings.

\begin{definition}
Let $A$ be a PI-algebra over a field $\k$. Its \emph{PI-exponent}, denoted by $\operatorname{exp}(A)$, is the integer
\[
\operatorname{exp}(A)=\lim_{n\to\infty}\sqrt[n]{c_n(A)},
\]
where $\{c_n(A)\}_{n\in\mathbb N}$ is the codimension sequence of $A$, defined by
\[
c_n(A)=
\dim_{\k}
\frac{P_n(\k)}
{P_n(\k)\cap T(A)},
\]
and $P_n(\k)$ denotes the space of multilinear polynomials in the variables $x_1,\ldots,x_n$ in $\k\langle X\rangle$.
\end{definition}

\begin{remark}
The existence of the limit defining $\operatorname{exp}(A)$ and the fact that it is always an integer are highly nontrivial results in the theory of PI-algebras. A detailed exposition, including proofs, can be found in \cite{GZ}.
\end{remark}

\begin{corollary}\label{PI-exp}
If $U$ is a Galois ring satisfying the hypotheses of Theorem~\ref{PI}, then
\[
\operatorname{exp}(U)=n^2,
\]
where $n=\operatorname{PI-deg}(U)$.
\end{corollary}

\begin{proof}
Since $T(U)=T(\mathcal{K})$ by Theorem~\ref{PI-0}, and $\mathcal{K}$ is a central simple algebra by Theorem~\ref{PI}, the result follows from \cite[Theorem 6.6.1]{GZ}.
\end{proof}

Our next theorem is an analogue for PI Galois rings of certain famous results in commutative algebra.

\begin{theorem}\label{integrality}
In the setting of Theorem \ref{PI}, suppose that $U$ is a Galois order. Then it is integral over its center $Z$. Hence
\begin{enumerate}[{\rm (i)}]
    \item (Lying over) If $\p \in \Spec \, Z$, there is $P \in \Spec \, U$ such that $P \cap Z=\p$.
    \item (Going up) If $\p_1, \p_2 \in \Spec \, Z$, $\p_1 \subset p_2$, then if $P \in \Spec \, U$ is such that $P \cap Z = \p_1$, there is $Q \in \Spec \, U$ containing $P$ such that $Q \cap Z = \p_2$.
    \item (Incomparability) If $P, Q \in \Spec \, U$ and $P \subset Q$, then $P \cap Z \subset Q \cap Z$.
    \item (Going down) If $\p_1, \p_2 \in \Spec \, Z$, $\p_1 \supset p_2$, then if $P \in \Spec \, U$ is such that $P \cap Z = \p_1$, there is $Q \in \Spec \, U$ contained in $P$ such that $Q \cap Z = \p_2$.
\end{enumerate}
\end{theorem}
\begin{proof}
    By Theorem \ref{center}(iii), $Z=\Gamma^\M$, and hence by Noether's Theorem in invariant theory, $\Gamma$ is an affine $Z$-algebra; by the very definition of a Galois ring, then $U$ is an affine $Z$-algebra, and then by \cite[Lemma 13.8.4(ii)]{McConnell}, $U$ is an integral $Z$-algebra. $U$, as we saw, is a prime PI-ring, and $\Gamma$, and hence $Z=\Gamma^\M$, is integrally closed. So, the rest of our Theorem follows from \cite[Theorem 13,8.14]{McConnell}.
\end{proof}



To finish this section, we will now consider an example of a PI Galois ring: the nilHecke algebra $\mathcal{H}(\mathfrak{h}, \mathcal{W})$  \cite{G18}, in the terminology of \cite{Webster} and \cite{Weekes}. We will follow \cite[Section 7.1 ]{G18}.
Let $\Sigma \subset \mathcal{R} \subset \mathfrak{h}^*$ be a reduced root system with simple roots in $\Sigma$, and let $\mathcal{W}$ be the Weyl group. For each $\alpha \in \mathcal{R}$ and the corresponding reflection $s_\alpha$ (viewed as an element of $\C \mathcal{W}$), we define $\theta_\alpha$ of $\C \mathcal{W} \ltimes \C(\mathfrak{h}^*)$ as follows: $\theta_\alpha:=\frac{1}{\alpha}(s_\alpha-1)$.  This can be extended to an injection $\C \mathcal{W} \into \C \mathcal{W} \ltimes \C(\mathfrak{h}^*)$. The algebra $\mathcal{H}(\mathfrak{h}, \mathcal{W})$ is the subalgebra of $\C \mathcal{W} \ltimes \C(\mathfrak{h}^*)$ generated by $\C[\mathfrak{h}^*]$ with basis $\theta_w, w \in \mathcal{W}$.

\begin{proposition}\label{nil-Hecke}
    $\mathcal{H}(\mathfrak{h}, \mathcal{W})$ is a Galois $\C[\mathfrak{h}^*]$-ring in $\C(\mathfrak{h}^*) * \C \mathcal{W}$. It is  a principal Galois order and a prime PI-algebra. 
\end{proposition}
\begin{proof}
    The first claim is clear from \cite[Propostion 2.9]{Hartwig} or \cite[Proposition 4.1(1)]{FO}. The second is \cite[Theorem 7.1.4]{G18} since the group is finite. The algebra is a PI and a prime ring due to Theorem \ref{localization4} and Proposition \ref{proposition-Fisher-Montgomery}.
\end{proof}

Since $\mathcal{H}(\mathfrak{h}, \mathcal{W})$ is free of rank $|\mathcal{W}|$ over its center $\C[h^*]^\mathcal{W}$ (cf. Theorem \ref{center}), applying the Chevalley-Shephard-Todd Theorem and \cite[Lemma 7.1.5, formula (7.1.6)]{G18} we get:

\begin{corollary}\label{nil-Hecke-2} 

The T-ideal of polynomial identities of $\mathcal{H}(\mathfrak{h}, \mathcal{W})$  is the same as the T-ideal of polynomial identities of $M_{|\mathcal{W}|}(\C)$. Its PI-degree is $|\mathcal{W}|$, and PI-exponent $|\mathcal{W}|^2$.
\end{corollary}
\begin{proof}
    This follows from the previous discussion, Proposition \ref{nil-Hecke} and Theorem \ref{PI}. The claim about the PI-exponent follows from Corollary \ref{PI-exp}.
\end{proof}

\section{Ring-theoretic dimensions of Galois rings}

In this section, except when discussing the noncommutative Krull dimension in the sense of Gabriel and Rentschler, we assume that all rings are algebras over a base field $\k$, which is otherwise arbitrary. The five dimensions considered in this section are the Gelfand--Kirillov dimension $\GK$, the classical Krull dimension $\cl$, the Krull dimension in the sense of Gabriel and Rentschler $\KK(\cdot)$, the Gelfand--Kirillov transcendence degree $\Tdeg$, and the lower transcendence degree $\LD$.

We begin by recalling the classical definition from commutative algebra, which extends naturally to possibly noncommutative rings.

\begin{definition}
A (not necessarily commutative) ring $R$ has classical Krull dimension $m$, written $\cl(R)=m$, if there exists a chain of prime ideals
\[
P_0 \subsetneq P_1 \subsetneq \cdots \subsetneq P_m
\]
of length $m$, and no such chain of greater length exists. If there is no upper bound on the lengths of such chains, we write $\cl(R)=\infty$.
\end{definition}

\subsection{Growth dimensions}

We now turn to the growth-theoretic invariants. For simplicity, we assume throughout that all algebras are affine.

\begin{definition}\cite{Gelfand}
Let $A$ be an affine algebra and let $V$ be a finite-dimensional vector space that generates $A$ as an algebra and contains $1$. Such subspaces are called \emph{frames}.

Let $d_V(n)=\dim V^n$. The Gelfand--Kirillov dimension of $A$ is defined by
\[
\GK A=\limsup_{n\to\infty}\log_n(d_V(n))
\footnote{$\log_n d_V(n)$ is an abbreviation for $\log d_V(n)/\log n$.}
\]
\end{definition}

It is well known that this definition is independent of the choice of frame (see \cite[Lemma 1.1]{KL}).

Before proceeding, we recall the following fundamental result.

\begin{theorem}\label{GK-comm}
If $A$ is a commutative $\k$-algebra, then $\GK A=\cl A$. If, in addition, $A$ is a domain, then both quantities are equal to the transcendence degree of $\Frac A$ over $\k$.
\end{theorem}

\begin{proof}
See \cite[Theorem 4.5]{KL}.
\end{proof}

The Gelfand--Kirillov transcendence degree and the lower transcendence degree are noncommutative analogues of the usual transcendence degree for division algebras. There is a plethora of such notions, each with its own advantages and drawbacks (see the discussions in \cite{Zhang,Zhang2,FSS,YZ}).

\begin{definition}\cite{Gelfand}
Let $A$ be an affine algebra. The Gelfand--Kirillov transcendence degree of $A$ is defined by
\[
\Tdeg A=
\inf_b
\limsup_{n\to\infty}
\log_n
\dim(\k+bV)^n,
\]
where $V$ is a frame of $A$ and $b$ ranges over all regular elements of $A$.
\end{definition}

Again, this definition is independent of the choice of frame (cf.~\cite{Zhang}). Although it is primarily used as an invariant of division algebras, it is defined for all associative algebras. If $A$ is a commutative field extension of $\k$, then $\Tdeg A$ coincides with the usual transcendence degree \cite[Proposition 2.2]{Zhang}, making it a natural candidate for a noncommutative analogue of transcendence degree.

Moreover, $\Tdeg$ is better suited than $\GK$ as an invariant of division algebras. Indeed, Makar-Limanov showed that the first Weyl field over $\mathbb{C}$ contains a free algebra on two generators and therefore has infinite Gelfand--Kirillov dimension, whereas its Gelfand--Kirillov transcendence degree is equal to $2$ (see \cite[Chapter 8]{KL}).

The lower transcendence degree, denoted $\LD$, was introduced in \cite{Zhang2}. It enjoys better ring-theoretic properties than the Gelfand--Kirillov transcendence degree and serves as an analogue of transcendence degree in noncommutative projective geometry (see \cite[Section 9]{Zhang2}). It is even conceivable that these two invariants coincide for division algebras, although they are known to differ for more general algebras. Whether they are always equal for division algebras remains an important open problem.

Let $A$ be an associative algebra and let $V$ be a frame of $A$. If, for every frame $V$, there exists a finite-dimensional vector subspace $W \subseteq A$ such that
\[
\dim VW=\dim W,
\]
then we define $\LD A=0$.

Otherwise, there exists a frame $V$ of $A$ such that, for every finite-dimensional subspace $W$,
\[
\dim VW \geq \dim W+1.
\]
We say that $VDI(A)_d$ holds for some real number $d>0$ (that is, $A$ \emph{satisfies the volume difference inequality of degree} $d$) if there exists a constant $c\in\mathbb{R}_{>0}$ such that, for every finite-dimensional subspace $W$,
\[
\dim VW \geq (\dim W+c\dim W)^{d(d-1)}.
\]

If instead
\[
\dim VW\geq \dim W+c\dim W,
\]
we say that $VDI(A)_\infty$ holds for $A$.

\begin{definition}
The lower transcendence degree of an affine algebra $A$ is defined as follows. If $\LD A\neq 0$, then
\[
\LD A=
\sup_V
\left\{
d \;\middle|\;
VDI(A)_d \text{ holds}
\right\},
\]
where $V$ ranges over all frames of $A$.
\end{definition}

An important notion introduced in \cite{Zhang2} is that of \emph{LD-stability}. For every affine algebra $A$,
\[
\LD A\leq \Tdeg A\leq \GK A.
\]
When equality holds throughout, we say that $A$ is \emph{LD-stable}. Prime PI algebras are LD-stable \cite{Zhang}. Hence, in particular, if $F$ is a finitely generated field extension of $\k$, then
\[
\LD F=\Tdeg F=\operatorname{tr.deg}_{\k}F,
\]
which is one of the desirable properties of the lower transcendence degree. Moreover, if $A$ is a prime Goldie ring, then
\[
\LD A=\LD \Q(A)=\Tdeg \Q(A),
\]
which considerably simplifies the computation of $\LD$ and $\Tdeg$ for division algebras, a task that is typically quite difficult.

We recall the following results.

\begin{proposition}\label{prop-zhang}
\begin{enumerate}[{\rm(i)}]
\item \emph{\cite[Proposition 2.1]{Zhang2}} Let $A$ be an affine algebra and let $S$ be a left or right denominator set consisting of regular elements. Then
\[
\LD A=\LD A_S.
\]

\item \emph{\cite[Theorem 0.3]{Zhang2}} Let $B\subseteq A$ be prime Goldie rings. Then
\[
\LD B\leq \LD A.
\]
Moreover,
\[
\LD A=\LD \Q(A).
\]

\item \emph{\cite[Proposition 0.4(2)]{Zhang2}} If $A$ is a prime PI-algebra, then
\[
\LD A=\Tdeg A=\GK A.
\]
\end{enumerate}
\end{proposition}

In \cite[Section 5]{FSS}, two of the present authors studied the lower transcendence degree of Galois rings under the assumptions that the base field is algebraically closed and of characteristic zero. We now remove both of these restrictions.

We begin with the following result, which is of independent interest.

\begin{proposition}\label{new-LD}
Let $A$ be an LD-stable algebra and let $B$ be a somewhat commutative algebra (cf.~\cite[Section 8.6]{McConnell}). Then $A\otimes B$ is LD-stable. Moreover,
\[
\GK(A\otimes B)=\GK A+\GK B.
\]
\end{proposition}

\begin{proof}
By \cite[Theorem 4.3(4)]{Zhang2}, the algebra $B$ is LD-stable. Furthermore, \cite[p.~28]{KL} shows that
\[
\GK(A\otimes B)=\GK A+\GK B.
\]
By \cite[Corollary 3.5 and Corollary 4.5]{Zhang2},
\[
\LD A+\LD(\gr B)
\leq
\LD(A\otimes B)
\leq
\LD A+\LD B.
\]
Since $\gr B$ is commutative,
\[
\LD(\gr B)=\GK(\gr B),
\]
and by \cite[Proposition 6.6]{KL},
\[
\GK(\gr B)=\GK B.
\]
Finally, since $A$ is LD-stable,
\[
\LD A=\GK A.
\]
Combining these equalities yields
\[
\LD(A\otimes B)
=
\GK(A\otimes B)
=
\GK A+\GK B,
\]
so $A\otimes B$ is LD-stable.
\end{proof}

\begin{corollary}\label{useful-degrees}
The Gelfand--Kirillov dimension and the lower transcendence degree of
$S_m^n(\k)$ and $Q_m^n(\k)$ are both equal to $m+n$.
\end{corollary}

\begin{proof}
By \cite[Theorem 0.5]{Zhang2} and the proposition above, the algebras
\[
W_m(\k)\otimes \k[z_1,\ldots,z_{n-m}]
\]
and
\[
\k_q[x,y]^{\otimes m}\otimes
\k[z_1,\ldots,z_{n-m}]
\]
are LD-stable and have Gelfand--Kirillov dimension $m+n$ (see, for example, \cite{McConnell,Brown}). Since
\[
\Frac S_m^n(\k)
=
\Frac\!\left(
W_m(\k)\otimes
\k[z_1,\ldots,z_{n-m}]
\right)
\]
and
\[
\Frac Q_m^n(\k)
=
\Frac\!\left(
\k_q[x,y]^{\otimes m}
\otimes
\k[z_1,\ldots,z_{n-m}]
\right),
\]
the result follows from Proposition~\ref{prop-zhang}(ii).
\end{proof}

\begin{remark}
As the last named author learned from a private communication with Ken Goodearl, it remains a difficult open problem to compute the Gelfand--Kirillov dimension of these algebras. Indeed, Gelfand--Kirillov dimension behaves poorly under localization in general, and the few known positive results do not apply to the present situation (see \cite[Chapters 4 and 12]{KL}).
\end{remark}

\begin{definition}
Let $A$ be an algebra and let $\alpha$ be an automorphism of $A$. We say that $\alpha$ is \emph{locally algebraic} if every element $r\in A$ is contained in an $\alpha$-stable finite-dimensional vector subspace of $A$.
\end{definition}

The following theorem is a simplified version of \cite[Theorem 6.1]{FO}.

\begin{theorem}\label{loc-alg}
Let $U$ be a Galois $\Gamma$-ring over an algebraically closed field of characteristic zero. Assume that $\mathcal{M}$ is a group generated by locally algebraic automorphisms $\alpha_1,\ldots,\alpha_j$. Then
\[
\GK U\geq \GK \Gamma+\operatorname{growth}\mathcal{M}.
\]
\end{theorem}

\begin{proof}
Let $V$ be a finite-dimensional vector subspace of $\Lambda$ with basis
$v_1,\ldots,v_s$. For each $m=1,\ldots,s$ and
$n=1,\ldots,j$, let $W_{mn}$ be a finite-dimensional vector subspace of
$\Lambda$ containing $v_m$ and stable under the action of $\alpha_n$.
Set
\[
W_m=\sum_{n=1}^{j}W_{mn},
\qquad
W=\sum_{m=1}^{s}W_m.
\]
Then $\mathcal{M}V\subseteq W$. Therefore, the hypotheses of
\cite[Theorem 6.1]{FO} are satisfied, and the result follows.
\end{proof}

\begin{corollary}
The alternating analogue of $U(\mathfrak{gl}_n)$,
$\mathfrak{A}(\mathfrak{gl}_n)$
(\cite[Definition 2.1]{Jauch}),
has Gelfand--Kirillov dimension at least $n^2$.
If $n=2$, then equality holds. Moreover,
$\mathfrak{A}(\mathfrak{gl}_2)$ is LD-stable.
\end{corollary}

\begin{proof}
The algebra $\Gamma$ has Gelfand--Kirillov dimension
$
\frac{n(n+1)}{2}
$
by \cite[Section 2.2, formula (5)]{Jauch}, and it is embedded into the invariant subalgebra of
$S_{n^2}^{n(n-1)/2}(\mathbb{C})$. Hence the first assertion follows from Theorem~\ref{loc-alg}.

If $n=2$, then $\mathfrak{A}(\mathfrak{gl}_2)$ is a finite free module over
$U(\mathfrak{gl}_2)$ by \cite[Lemma 3.1]{Jauch}. Therefore,
\cite[Proposition 8.2.9(ii)]{McConnell} and
\cite[Proposition 3.1]{Zhang2} imply that
\[
\GK \mathfrak{A}(\mathfrak{gl}_2)
=
\LD \mathfrak{A}(\mathfrak{gl}_2)
=
\GK U(\mathfrak{gl}_2)
=
4.
\]
This proves the claim.
\end{proof}

\begin{definition}\cite{Zhang}
An affine algebra $A$ for which the classical ring of fractions $\Q(A)$ exists is called \emph{Tdeg-stable} if
\[
\Tdeg \Q(A)=\Tdeg A=\GK A.
\]
\end{definition}

\begin{proposition}
Let $U$ be a $\Gamma$-ring. If $\mathcal{M}$ is a finite group, then
\[
\GK U=\GK \Gamma,
\]
and $U$ is a PI ring. Moreover, if $|G|$ and $|\mathcal{M}|$ are invertible in
\[
\mathcal{L}=L*\mathcal{M},
\]
then $U$ is a semiprime PI Goldie ring. In this case, $U$ is also Tdeg-stable.
\end{proposition}

\begin{proof}
The ring $\mathcal{L}$ is clearly a PI ring, since $L*\mathcal{M}$ is finite-dimensional over $L$. Therefore, $U$ is also a PI ring. We have a chain of inclusions
\[
\Gamma \subseteq U \subseteq \mathcal{K} \subseteq L*\mathcal{M}.
\]
By \cite[Proposition 8.2.9(ii)]{McConnell},
\[
\GK(L*\mathcal{M})=\GK L.
\]
Furthermore, by \cite[Proposition 2.2]{Zhang},
\[
\GK L=\operatorname{tr.deg}L
=\operatorname{tr.deg}K
=\cl\Gamma,
\]
while Theorem~\ref{GK-comm} yields
\[
\cl\Gamma=\GK \, \Gamma.
\]
Hence
\[
\GK\Gamma
\leq
\GK U
\leq
\GK(L*\mathcal{M})
=
\GK\Gamma,
\]
and therefore
$
\GK U=\GK\Gamma.
$

Since $\mathcal{M}$ is finite and $|\mathcal{M}|$ is invertible, it follows from \cite[Corollary 0.2]{Montgomery} that $\mathcal{L}$ is semisimple Artinian. Moreover, \cite[Lemma 1.13]{Montgomery} implies that $\mathcal{K}$ is semisimple Artinian. Hence, by Theorem~\ref{localization4}, $U$ is a semiprime Goldie ring. Finally, \cite[Theorem 1.1(2)]{Zhang} implies that $U$ is Tdeg-stable.
\end{proof}

\begin{corollary}\label{coincide-1}
If $U$ is an affine prime PI Galois $\Gamma$-ring with center $Z$, then
\[
\GK \, U=\GK \, \Gamma=\GK \, Z.
\]
In particular, if
$
F=\Frac Z,
$
then $K$ is an algebraic extension of $F$.
\end{corollary}

\begin{proof}
The statement follows from the proposition above, Theorem~\ref{PI}, and \cite[Theorems 1.13.7 and 1.13.9]{GZ}.
\end{proof}

As an illustration of the applicability of our methods, we briefly discuss the noncommutative deformations of Kleinian singularities introduced in \cite{C-B-H}. In type $A$, Hodges \cite{Hodges} showed that every such deformation is a generalized Weyl algebra of rank one and hence a Galois order. The results of \cite{FSS} imply that these algebras are LD-stable domains of Gelfand--Kirillov dimension $2$. Hartwig \cite{Hartwig2} proved that the deformations of type $D$, which we denote by $D(q)$ for a polynomial $q\in\mathbb{C}[t]$, are principal Galois orders.

\begin{theorem}\label{D-deformation}
The algebra $D(q)$ admits a natural finite-dimensional filtration
$\mathcal{F}$ that makes it a somewhat commutative algebra. Hence,
$D(q)$ is LD-stable, and
\[
\LD D(q)=\Tdeg D(q)=\GK D(q)=2.
\]
\end{theorem}

\begin{proof}
The filtration constructed in \cite[Remark 2.3]{Hartwig2}, together with
\cite[Proposition 2.5(i)]{Hartwig2}, implies that $D(q)$ is a somewhat
commutative algebra. Hence, $D(q)$ is LD-stable by
\cite[Theorem 4.3(4)]{Zhang2}. The associated graded algebra is
$\mathbb{C}[x,y]^{D_n}$, where $D_n$ is the binary dihedral group. By
\cite[Proposition 6.6]{KL},
\[
\GK D(q)
=
\GK\bigl(\mathbb{C}[x,y]^{D_n}\bigr)
=
2,
\]
since the latter is an affine commutative domain of Krull dimension~$2$.
\end{proof}

We conclude this section by computing the Gelfand--Kirillov dimension of
the Galois ring $U_q(\mathfrak{gl}_n)$ (\cite{Hartwig}) and showing that
it is LD-stable.

\begin{theorem}\label{GK-quantum-group}
If $\k$ is an algebraically closed field of characteristic zero and
$q$ is not a root of unity, then $U_q(\mathfrak{gl}_n)$ is LD-stable and
has Gelfand--Kirillov dimension $n^2$.
\end{theorem}

\begin{proof}
As observed in \cite{FH2}, the algebra $U_q(\mathfrak{gl}_n)$ admits a
multi-filtration by $\mathbb{N}^s$, for some $s\in\mathbb{N}$, whose
associated graded algebra is a localization of a quantum affine space
(\cite{Brown}) in $n^2$ indeterminates at local normal commuting
indeterminates, in the terminology of \cite{LMO}. Combining the
computation of the Gelfand--Kirillov dimension of quantum affine spaces
from \cite{Brown} with the localization theorem
\cite[Theorem 2]{LMO}, we obtain
\[
\GK\bigl(\gr U_q(\mathfrak{gl}_n)\bigr)=n^2.
\]
It also follows from the results of \cite{Zhang2} that
$\gr U_q(\mathfrak{gl}_n)$ is LD-stable. Since $\mathbb{N}^s$ is an
ordered semigroup, \cite[Theorem 4.3(4)]{Zhang2} implies that
$U_q(\mathfrak{gl}_n)$ is LD-stable and has
Gelfand--Kirillov dimension $n^2$.
\end{proof}

\subsection{Krull dimension}

We now introduce the Krull dimension in the sense of
Gabriel--Rentschler (cf.~\cite[Chapter 15]{GW} and
\cite[Chapter 6]{McConnell}).

\begin{definition}
Let $R$ be a ring and $M$ an $R$-module. Using transfinite induction, we
define the (left) Krull dimension in the sense of
Gabriel--Rentschler, denoted by $\KK(M)$.

Set
\[
\KK_{-1}^R=\{0\}.
\]
Let $\alpha\geq 0$ be an ordinal, and suppose that
$\KK_\beta^R$ has been defined for every ordinal
$\beta<\alpha$. Then an $R$-module $M$ belongs to
$\KK_\alpha^R$ if, for every countable descending chain of submodules
\[
M_0\supset M_1\supset M_2\supset\cdots,
\]
the quotient modules
$
M_i/M_{i+1}
$
belong to $\KK_\gamma^R$ for some $\gamma<\alpha$, for all but finitely
many $i$.

The Krull dimension of a ring $R$, denoted by $\KK(R)$, is the Krull
dimension of $R$ regarded as a left $R$-module.
\end{definition}

For instance, a module has Krull dimension~0 if and only if it is nonzero and Artinian. Unlike the other ring-theoretic dimensions, the Krull dimension is not always defined. For example, it is undefined for modules containing an infinite direct sum of a nonzero submodule \cite[Exercise 15.C]{GW}. On the other hand, $\KK(R)$ is always defined for a Noetherian ring $R$ \cite[Lemma 15.3]{GW}.

In general, it is an open problem whether the left and right Krull dimensions of a Noetherian ring always coincide.

\begin{theorem}\label{coincidence-of-dimension-2}
Let $A$ be a Noetherian affine PI algebra over a field $\k$. If $\GK A$ is finite, then
\[
\cl(A)=\KK(A)=\GK(A).
\]
\end{theorem}

\begin{proof}
The statement follows from a combination of
\cite[Corollary 3.5.52]{Rowen},
\cite[Proposition 6.1.42]{Rowen2}, and
\cite[Corollary 10.16(b)]{KL}.
\end{proof}

\begin{corollary}\label{prime-PI-Galois}
Let $U$ be an affine prime Noetherian PI algebra over a field $\k$ that is a Galois $\Gamma$-ring with center $Z$. Then
\[
\cl(U)=\KK(U)=\GK(U)=\GK(\Gamma)=\GK(Z).
\]
In particular, if
$
F=\Frac Z,
$
then $K$ is an algebraic extension of $F$.
\end{corollary}

\begin{proof}
The result follows from Theorem~\ref{coincidence-of-dimension-2} and Corollary~\ref{coincide-1}.
\end{proof}

The following theorem summarizes the basic properties of the Krull dimension.

\begin{theorem}\label{Main-Krull}\phantom{X}
\begin{enumerate}[{\rm(i)}]

\item
Let $R$ be a ring and let $X$ be a denominator set. Then
\[
\KK(R_X)\leq\KK(R).
\]

\item
Let $R\subseteq S$ be Noetherian rings such that $S$ is faithfully flat as a left $R$-module. In particular, this holds if $S$ is a free $R$-module. Then
\[
\KK(R)\leq\KK(S).
\]

\item
Let $R$ be a ring equipped with an $\mathbb{N}$-filtration $\mathcal{F}$. Then
\[
\KK(R)\leq\KK(\gr_{\mathcal{F}}R).
\]

\item
Let $R$ be a Noetherian ring and let $\theta$ be an automorphism of $R$. Then
\[
\KK(R[x;\theta])=\KK(R)+1.
\]

\item
The Krull dimension is a Morita invariant.

\item
Let $R\subseteq S$ be Noetherian rings such that $S$ is a finitely generated free $R$-module. Then
\[
\KK(R)=\KK(S).
\]

\item
Let $G$ be a finite group of automorphisms of a Noetherian ring $R$ such that $|G|^{-1}\in R$. Then
\[
\KK(R^G)\leq\KK(R).
\]

\end{enumerate}
\end{theorem}

\begin{proof}
Item~(i) is proved in \cite[Exercise 15S]{GW}, item~(ii) in
\cite[Exercise 15U]{GW}, item~(iii) in
\cite[Lemma 6.5.6]{McConnell}, item~(iv) in
\cite[Theorem 15.19]{GW}, item~(v) in
\cite[Proposition 6.5.1(ii)]{McConnell}, item~(vi) in
\cite[Corollary 6.5.3]{McConnell}, and item~(vii) in
\cite[Theorem 7.8.8(i)]{McConnell}.
\end{proof}

We now apply these results to Galois orders.

\begin{lemma}\label{krull-invariants}
Let $R$ be a simple Noetherian ring, and let $G$ be a finite group acting faithfully on $R$ by outer automorphisms. Assume that $|G|^{-1}\in R$. Then
\[
\KK(R^G)=\KK(R).
\]
\end{lemma}

\begin{proof}
By Theorem~\ref{invariant-theory}, the rings $R^G$ and $R*G$ are Morita equivalent. Since $G$ is finite, the skew group ring $R*G$ is Noetherian. As $\KK(\cdot)$ is a Morita invariant,
\[
\KK(R^G)=\KK(R*G).
\]
On the other hand, Theorem~\ref{Main-Krull}(vi) yields
\[
\KK(R*G)=\KK(R).
\]
The result follows.
\end{proof}

\begin{proposition}\label{Krull-skew-group}
Let $U$ be a Galois ring in $\mathcal{L}=L*\mathbb{Z}^n$. Then
\[
\KK(\mathcal{L})=\KK(\mathcal{K})=n, \,\,\, \KK(U)\geq n.
\]
\end{proposition}

\begin{proof}
By Theorem~\ref{U-is-prime-Goldie}, $\mathcal{L}$ is a simple Noetherian ring. By \cite[Proposition 13]{Park},
\[
\KK(\mathcal{L})=h(\mathbb{Z}^n),
\]
where $h(\mathbb{Z}^n)$ denotes the Hirsch number of $\mathbb{Z}^n$ (cf.~\cite[Definition 8.2.9]{Rowen2}). By \cite[Proposition 8.2.10]{Rowen2},
\[
h(\mathbb{Z}^n)=n.
\]
Hence Lemma~\ref{krull-invariants} implies that
$
\KK(\mathcal{K})=n.
$
Finally, Theorem~\ref{Main-Krull}(i), together with Theorem~\ref{localization2}, gives
$
\KK(U)\geq n.
$
\end{proof}




\section{Affine Iwahori--Hecke algebras and DAHAs}

Throughout this section, we assume that $\k=\mathbb{C}$. We now discuss affine Iwahori--Hecke algebras and double affine Hecke algebras (DAHAs), following the approach of \cite{GKV}. Their connection with Galois orders was pointed out in \cite[Example 3.8 and Theorem 3.9]{Hartwig2}. We repeat the discussion here in greater detail in order to illustrate our methods.

Let $A=(a_{ij})_{n\times n}$ be a symmetrizable generalized Cartan matrix. Associated with $A$ is a root datum consisting of a free abelian group $X$ of rank $2n-\operatorname{rk}(A)$, its dual lattice
\[
X^\vee=\operatorname{Hom}(X,\mathbb{Z}),
\]
a perfect pairing
\[
\langle\, , \,\rangle:X^\vee\times X\longrightarrow\mathbb{Z},
\]
together with sets of simple roots
\[
\alpha_1,\ldots,\alpha_n\in X
\]
and simple coroots
\[
\alpha_1^\vee,\ldots,\alpha_n^\vee\in X^\vee.
\]

The vector space
$
\mathbb{C}\otimes_{\mathbb{Z}}X,
$
together with the simple roots and coroots, is a realization of the Kac--Moody Lie algebra $\mathfrak{g}(A)$. The Cartan subalgebra is
$
\mathfrak{h}=\mathbb{C}\otimes_{\mathbb{Z}}X^\vee,
$
and the corresponding maximal torus of the Kac--Moody group is
$
\mathbb{C}^*\otimes_{\mathbb{Z}}X^\vee.
$


Let $R$ be the root system of $\mathfrak{g}(A)$ and $W$ its Weyl group. Set
$
T_{\mathcal{C}}:=\mathcal{C}\otimes_{\mathbb{Z}}X^\vee,
$
where $\mathcal{C}$ is a normal algebraic affine curve (for example, a smooth affine curve). In \cite[Definition 1.3 and Theorem 1.4]{GKV}, the authors introduced the algebra $\mathsf{H}_q(\mathcal{C},A)$, where $q\in\mathbb{C}$ and $q\neq\pm1$, whose definition we now recall.

\begin{definition}\label{def-Hecke}
Consider the skew group ring $\mathbb{C}(T_{\mathcal{C}})*W$. The algebra $\mathsf{H}_q(\mathcal{C},A)$ is the subalgebra spanned by the elements
\[
\sum_{w\in W}f_w\,w
\]
such that:

\begin{enumerate}[{\rm(i)}]

\item
Each function $f_w$ has no singularities other than simple poles along the divisors $T_{\mathcal{C}}^\alpha$, for finitely many $\alpha\in R_+^{\mathrm{re}}$.

\item
Given $\alpha\in R_+^{\mathrm{re}}$ and $x\in\mathcal{C}$, let $T_{\alpha,x}$ denote the subvariety of $T_{\mathcal{C}}$ consisting of those points $t\in T_{\mathcal{C}}$ such that, if $t=c\otimes\nu$, where $c\in\mathcal{C}$ and $\nu\in X^\vee$, then
\[
c^{\langle\alpha,\nu\rangle}=x.
\]

\item
The element
\[
\sum_{w\in W}f_w\,w
\]
belongs to $\mathsf{H}_q(\mathcal{C},A)$ if, for every $w\in W$ and every $\alpha\in R_+^{\mathrm{re}}$,
\[
\operatorname{Res}_{T_{\alpha,1}}(f_w)
+
\operatorname{Res}_{T_{\alpha,1}}(f_{s_\alpha w})
=
0.
\]

\item
The function $f_w$ vanishes on $T_{\alpha,q^{-2}}$ whenever $\alpha\in R_+^{\mathrm{re}}$ and
$
w^{-1}(\alpha)\in R_-.
$

\end{enumerate}
\end{definition}

The following theorem has the same proof as
\cite[Example 3.8 and Theorem 3.9]{Hartwig2}.

\begin{theorem}\label{Hecke-1}
The algebra $\mathsf{H}_q(\mathcal{C},A)$ is a principal
$\mathbb{C}[T_{\mathcal C}]$-Galois order in
$\mathbb{C}(T_{\mathcal C})*W$.
\end{theorem}

Suppose that $A$ is of finite type.

\begin{theorem}\label{Hecke-2}
The algebra $\mathsf{H}_q(\mathcal{C},A)$ is a Noetherian prime ring with
center $\mathbb{C}[T_{\mathcal C}]^W$.
\end{theorem}

\begin{proof}
By Theorem~\ref{U-is-prime-Goldie},
$\mathsf{H}_q(\mathcal{C},A)$ is a prime Goldie ring. Since it is
Noetherian by \cite[Theorem 2.5]{GKV}, it is a Noetherian prime ring.
Finally, Theorem~\ref{center} yields
\[
Z\bigl(\mathsf{H}_q(\mathcal{C},A)\bigr)
=
\mathbb{C}[T_{\mathcal C}]^W.
\]
\end{proof}

If $\mathcal{C}=\mathbb{C}^*$, respectively $\mathbb{C}$, then
$\mathsf{H}_q(\mathcal{C},A)$ is the affine, respectively degenerate
affine, Iwahori--Hecke algebra associated with the transpose Cartan
matrix $A^t$ (cf.~\cite[Section 5]{GKV}). Thus, we recover and extend
their basic structural properties using only the theory of Galois
orders.

The other case of interest is when $A$ is of untwisted affine type.
Let $A_0$ denote its finite part, with finite Weyl group $W$, root
system $R$, and coroot lattice $Q^\vee$
(cf.~\cite[Section 6]{GKV}). In this case,
\[
\mathfrak h
=
\mathfrak h_0
\oplus
\mathbb C c
\oplus
\mathbb C d,
\]
where $\mathfrak h_0$ is the Cartan subalgebra of the finite-dimensional
Lie algebra corresponding to $A_0$.

We have

\begin{theorem}\label{Hecke-4}
The double affine Hecke algebra associated with an untwisted affine
Cartan matrix $A$ is isomorphic to
$\mathsf{H}_q(\mathbb C^*,A)$. Consequently, it is a principal Galois
order in
\[
\mathbb C(\mathbb C c\oplus\mathfrak h_0)
*
(Q^\vee\rtimes W).
\]
Moreover, it is a prime ring with center $\mathbb C[c]$.
\end{theorem}

\begin{proof}  The first assertion is due to
\cite[Theorem 1.8]{GKV} and
\cite[Theorem 3.9]{Hartwig2}.
The computation of the center follows from
Theorem~\ref{center}; this description is also well known
(cf.~\cite[Section 6]{GKV}). The primeness follows from
Theorem~\ref{U-is-prime-Goldie}.
\end{proof}

\subsection{Gelfand--Kirillov Conjecture for the DAHA}

The Gelfand--Kirillov conjecture, now known to be false in general
(\cite{AOvdB}), originally asserted that, for an algebraic Lie algebra
$\mathfrak g$ over a field of characteristic zero, the skew field of
fractions $\Frac \, U(\mathfrak g)$ of its universal enveloping algebra $U(\mathfrak g)$ is
isomorphic to a Weyl skew field \cite{Gelfand}.\footnote{For
non-algebraic Lie algebras, the conjecture was already shown to fail in
\cite{Gelfand}.}

More generally, given an Ore domain $U$, one says that $U$ satisfies the
Gelfand--Kirillov conjecture (or \emph{Gelfand--Kirillov hypothesis}) if
$\Frac \, U$ is isomorphic to a Weyl skew field.

Several variations of this conjecture have been studied in different
settings (see \cite{SchwarzPan}), including symplectic reflection
algebras \cite{EG}, $W$-algebras \cite{FMO,Petukhov}, invariant rings of
differential operators \cite{AD,EFOS,FS,Tikaradze,SchwarzPan, Schwarz2}, and their
quantum analogues (see, for example, \cite{Brown,FH}).

Now let $\mathsf{H}_q(A)$ denote the double affine Hecke algebra
associated with an untwisted affine Cartan matrix $A$. In the notation
of Theorem~\ref{Hecke-4}, let $W$ be the finite Weyl group of the finite
part of $A$, and consider the idempotent
\[
e=1/|W|\sum_{w\in W}w.
\]
Define the spherical subalgebra
\[
\mathsf{U}_q(A):=e\,\mathsf{H}_q(A)\,e,
\]
whose identity element is $e$. Applying
\cite[Theorem 4.5]{Schwarz2} together with
Theorem~\ref{Hecke-4}, we obtain the following result.

\begin{proposition}\label{preliminar-GKC-DAHA}
The algebra
$
\mathsf{U}_q(A)
$
is a Galois order in
\[
\bigl(\mathbb{C}(\mathbb{C}c\oplus\mathfrak h_0)*Q^\vee\bigr)^W.
\]
Moreover, it is an Ore domain.
\end{proposition}

\begin{proof}
The first assertion follows from the preceding discussion. The fact that
$\mathsf{U}_q(A)$ is an Ore domain follows from
Theorem~\ref{U-is-Ore}.
\end{proof}

By Lemma~\ref{fundamental-lattices} and
\cite[Theorem 6.1]{FS}, we obtain the following theorem.

\begin{theorem}\label{GKC-DAHA}
The algebra $\mathsf{U}_q(A)$ satisfies the Gelfand--Kirillov
conjecture. More precisely, if $A$ is an $n\times n$ untwisted affine Cartan matrix,
then
\[
\Frac\mathsf{U}_q(A)
\simeq
\Frac W_{n-1}\bigl(\mathbb{C}(c)\bigr).
\]
\end{theorem}

\subsection{PI-theory of (double) affine Hecke algebras}

We now turn to polynomial identities.

\begin{theorem}\label{identities-Hecke}
Let $\mathsf{H}_q(A)$ be the double affine Hecke algebra associated with an untwisted affine Cartan matrix $A$. Then $\mathsf{H}_q(A)$ is not a PI algebra. If $A$ is of finite type, then the algebras $\mathsf{H}_q(\mathcal{C},A)$ are PI algebras. Moreover, their $T$-ideal of polynomial identities coincides with the $T$-ideal of polynomial identities of $M_{|W|}(\mathbb{C})$. Consequently, they have PI-degree $|W|$ and PI-exponent $|W|^2$. In particular, this applies to the affine Iwahori--Hecke algebra and its degenerate analogue.
\end{theorem}

\begin{proof}
The algebra $\mathsf{H}_q(A)$ is not a PI-algebra by Lemma~\ref{fundamental-lattices} and Theorem~\ref{usually-not-PI}.

Suppose now that $A$ is of finite type. By Theorems~\ref{Hecke-1} and \ref{PI-0}, the algebra $\mathsf{H}_q(\mathcal{C},A)$ is a PI-algebra, since the skew group algebra
$
\mathbb{C}(T_{\mathcal C})*W
$
is finite-dimensional over its center
$
\mathbb{C}(T_{\mathcal C})^W.
$
Hence its PI-degree is $|W|$, and its PI-exponent is $|W|^2$. Since $\mathsf{H}_q(\mathcal{C},A)$ is moreover a prime ring by Theorem~\ref{Hecke-2}, the conclusion follows from Theorem~\ref{PI}.
\end{proof}

The algebras $\mathsf{H}_q(\mathcal{C},A)$ are affine PI algebras. Hence, by the Razmyslov--Kemer--Braun Theorem (see \cite[Theorem 1.11.14]{GZ}), their Jacobson radical is nilpotent, and as the algebras are prime, it is zero:

\[
\operatorname{Jac}\bigl(\mathsf{H}_q(\mathcal{C},A)\bigr)=0.
\]
Therefore, we obtain the following corollary.

\begin{corollary}\label{affine-Hecke-semiprimitive}
The algebras $\mathsf{H}_q(\mathcal{C},A)$ are semiprimitive.
\end{corollary}

We also obtain a complete description of the various ring-theoretic dimensions of $\mathsf{H}_q(\mathcal{C},A)$.

\begin{theorem}\label{dimensions-affine-Hecke}
Let $A$ be an $n\times n$ Cartan matrix of finite type. Then
\[
\cl\bigl(\mathsf{H}_q(\mathcal{C},A)\bigr)
=
\KK\bigl(\mathsf{H}_q(\mathcal{C},A)\bigr)
=
\GK\bigl(\mathsf{H}_q(\mathcal{C},A)\bigr)
=
\LD\bigl(\mathsf{H}_q(\mathcal{C},A)\bigr)
=
\Tdeg\bigl(\mathsf{H}_q(\mathcal{C},A)\bigr)
=
n.
\]
\end{theorem}

\begin{proof}
The algebras $\mathsf{H}_q(\mathcal{C},A)$ are affine Noetherian prime PI algebras. The result therefore follows by combining Proposition~\ref{prop-zhang}, Corollary~\ref{coincide-1}, and Theorem~\ref{coincidence-of-dimension-2}.
\end{proof}

\subsection{The Nullstellensatz for affine Hecke algebras}

There is no universal agreement in the literature on what should be
considered the appropriate noncommutative analogue of the
Nullstellensatz; see \cite[Section 2.12, Supplement K]{Rowen}. We shall
use the following notions, following \cite{Rowen} and
\cite{McConnell}.

\begin{definition}\label{usual-Nullstellensatz}
Let $A$ be a $\k$-algebra. We say that $A$ has the \emph{radical property}
if every factor algebra of $A$ has nil Jacobson radical. We say that $A$
has the \emph{endomorphism property} over $\k$ if, for every simple
$A$-module $M$, the division algebra $\End_A(M)$ is algebraic over
$\k$. If $A$ satisfies both the radical property and the endomorphism
property over $\k$, then we say that $A$ satisfies the
\emph{Nullstellensatz} over $\k$.
\end{definition}

\begin{definition}\label{maximal-Nullstellensatz}
Let $A$ be a $C$-algebra, where $C$ is a commutative ring. We say that
$A$ satisfies the \emph{maximal Nullstellensatz} over $C$ if, for every
simple $A$-module $M$, the ideal
\[
P=\operatorname{Ann}_C(M)
\]
is maximal in $C$, and the algebra
\[
A/\operatorname{Ann}_A(M)
\]
satisfies the endomorphism property over $C/P$.
\end{definition}

\begin{definition}
A ring $R$ is called a \emph{Jacobson ring} if every prime ideal is the
intersection of primitive ideals.
\end{definition}

The following proposition will be used repeatedly.

\begin{proposition}\label{NJ}
Let $A$ be a countably generated Noetherian algebra over an uncountable
field $\k$. Then $A$ satisfies the Nullstellensatz and is a Jacobson
ring.
\end{proposition}

\begin{proof}
By \cite[Corollary 9.1.8]{McConnell}, $A$ satisfies the
Nullstellensatz over $\k$. Since $A$ is Noetherian, every prime factor
ring of $A$ is again Noetherian and hence Goldie. Therefore,
\cite[Lemma 9.1.2(ii)]{McConnell} implies that $A$ is a Jacobson ring.
\end{proof}

\begin{theorem}\label{Nulls-affine-Hecke}
Suppose that $A$ is of finite type. Then the algebras
$\mathsf{H}_q(\mathcal{C},A)$ are Jacobson rings and satisfy the maximal
Nullstellensatz over their center
$\mathbb{C}[T_{\mathcal C}]^W$. They also satisfy the Nullstellensatz
over $\mathbb{C}$.
\end{theorem}

\begin{proof}
By Noether's theorem on invariants,
$\mathbb{C}[T_{\mathcal C}]^W$ is a finitely generated
$\mathbb{C}$-algebra and hence a Jacobson ring by
\cite[Corollary 2.12.37]{Rowen}. Since
$\mathsf{H}_q(\mathcal{C},A)$ is affine over its center and is a
PI algebra, \cite[Theorem 6.3.3]{Rowen2} implies that it is a Jacobson
ring and satisfies the maximal Nullstellensatz.

Finally, $\mathsf{H}_q(\mathcal{C},A)$ is finitely generated over the
uncountable field $\mathbb{C}$, so Proposition~\ref{NJ} shows that it
also satisfies the Nullstellensatz over $\mathbb{C}$.
\end{proof}

\subsection{Integrality}

As an immediate consequence of Theorem~\ref{integrality}, we obtain the
following result.

\begin{theorem}\label{Hecke-integral}
The algebras $\mathsf{H}_q(\mathcal{C},A)$ are integral over their
center $\mathbb{C}[T_{\mathcal C}]^W$. Consequently, they satisfy the
Lying Over, Going Up, Incomparability, and Going Down properties with
respect to their center.
\end{theorem}

We conclude this section by observing that a very similar analysis,
leading to analogous results, can also be carried out for the nilHecke
algebras $\mathcal{H}(\mathfrak h,\mathcal W)$.

 \section{Spherical Coulomb branch algebras}

In this section $\k=\mathbb{C}$. Let $G$ be a reductive algebraic group and
$T$ a maximal torus, with Weyl group $W$. Let $F$ be another torus,
which in the physics literature is called the flavour torus. Consider an
extension of $G$ by $F$:
\[
0 \longrightarrow G \longrightarrow \bar{G} \longrightarrow F
\longrightarrow 0 .
\]

Let $\bar{T}$ be a maximal torus of $\bar{G}$ containing $F$, and let
$N$ be a finite-dimensional $\bar{G}$-representation. In \cite{BFN},
a delicate construction of an associative commutative product on an
equivariant Borel--Moore homology group of an infinite-dimensional
scheme, followed by taking the functor $\operatorname{Spec}$, gives an
affine normal Poisson variety $\mathcal{M}_C(G,N)$, called the
(flavour deformation of the) Coulomb branch. This construction also
provides a natural quantization (obtained by adding the loop parameter
$\mathbb{C}^*$ to the equivariant homology), called in \cite{lots} the
spherical Coulomb branch algebra. We denote it by
$\mathcal{A}_{\hbar}(G,N)$ (for an explanation of this terminology, see
\cite{Weekes}). Throughout this section we specialize $\hbar=1$ and
write simply $\mathcal{A}(G,N)$.

Denote by $\bar{\mathfrak t}$, $\mathfrak t$, and $\mathfrak f$ the
abelian Lie algebras of the algebraic tori $\bar{T}$, $T$, and $F$,
respectively. Let
\[
\bar{W}=X^\vee\rtimes W
\]
be the extended affine Weyl group. One of the main results of Webster
\cite[Proposition 4.2]{Webster} is the following.

\begin{theorem}\label{Coulomb-is-Galois-ring}
The algebra $\mathcal{A}(G,N)$ is a principal
$\mathbb{C}[\bar{\mathfrak t}]^W$-Galois order in
\[
(\mathbb{C}(\bar{\mathfrak t})*X^\vee)^W .
\]
\end{theorem}

We recall some basic structural properties of these algebras.

\begin{proposition}\label{Noetherian-domain}
The algebra $\mathcal{A}(G,N)$ is a Noetherian (and hence Ore) domain.
\end{proposition}

\begin{proof}
As noted in \cite[Section 1.2]{lots}, $\mathcal{A}(G,N)$ has an
$\mathbb{N}$-filtration such that the associated graded algebra is the
algebra of regular functions on the Coulomb branch
$\mathcal{M}_C(G,N)$. Since this algebra is Noetherian and a domain,
standard filtered--graded arguments imply that $\mathcal{A}(G,N)$ is
also Noetherian and a domain (see, e.g.,
\cite[Chapter~1]{McConnell}).
\end{proof}

We remark that the fact that $\mathcal{A}(G,N)$ is an Ore domain also
follows from Theorem~\ref{U-is-Ore} and Theorem~\ref{Coulomb-is-Galois-ring}.

The structure theory of Galois rings allows us to show that the
subalgebra $\mathbb{C}[\bar{\mathfrak t}]^W$ is a Harish--Chandra
subalgebra and a maximal commutative subalgebra of $\mathcal{A}(G,N)$.

\begin{proposition}\label{Coulomb-maximal-commutative}
Let $\mathcal{A}(G,N)$ be a spherical Coulomb branch algebra. Then
$\mathbb{C}[\bar{\mathfrak t}]^W$ is a Harish--Chandra subalgebra and a
maximal commutative subalgebra.
\end{proposition}

\begin{proof}
The claim follows from Theorem~\ref{HC-subalgebra}.
\end{proof}

As observed in \cite{lots}, by \cite[Theorem 4.1]{FO} (or by our
Theorem~\ref{center}), the center of $\mathcal{A}(G,N)$ is
$\mathbb{C}[\mathfrak f]$. In particular, when the flavour group is
trivial, the center consists only of scalars. Moreover, by the results
of \cite{BFN}, the spherical Coulomb branch algebra is free over its
Harish--Chandra subalgebra.

\subsection{Basic ring-theoretic properties}

In this subsection we show that spherical Coulomb branch algebras share
many properties with enveloping algebras of finite-dimensional complex
Lie algebras.

Our first result in this direction is that these algebras satisfy the
Nullstellensatz over $\mathbb{C}$ and are Jacobson rings.

\begin{theorem}\label{Null+Jacobson}
Every spherical Coulomb branch algebra $\mathcal{A}(G,N)$ satisfies the
Nullstellensatz over $\mathbb{C}$ and is a Jacobson ring.
\end{theorem}

\begin{proof}
The algebra $\mathcal{A}(G,N)$ is a Noetherian algebra finitely generated
over the uncountable field $\mathbb{C}$. Hence the result follows from
Proposition~\ref{NJ}.
\end{proof}

The following corollary generalizes a familiar property of enveloping
algebras (cf. \cite[Theorem 9.4.14]{Rowen2}).

\begin{corollary}\label{CB-semiprimitive}
The Jacobson radical of $\mathcal{A}(G,N)$ is $(0)$.
\end{corollary}

\begin{proof}
Since $\mathcal{A}(G,N)$ is a domain, $(0)$ is a prime ideal, and hence
an intersection of primitive ideals. Since
$\operatorname{Jac}(\mathcal{A}(G,N))$ is the intersection of all
primitive ideals, it follows that
\[
\operatorname{Jac}(\mathcal{A}(G,N))=(0).
\]
\end{proof}

\begin{corollary}
Let $M$ be an irreducible
$\mathcal{A}(G_1,N_1)$-module and let $L$ be an irreducible
$\mathcal{A}(G_2,N_2)$-module. Assume that the flavour torus is trivial.
Then
$
M\otimes L
$
is an irreducible $\mathcal{A}(G,N)$-module, where
\[
G=G_1\times G_2,\qquad N=N_1\oplus N_2 .
\]
\end{corollary}

\begin{proof}
It is shown in \cite[3(vii)(a)]{BFN} that
\[
\mathcal{A}(G_1,N_1)\otimes\mathcal{A}(G_2,N_2)
\simeq
\mathcal{A}(G,N).
\]
The result now follows from \cite[Theorem 1]{Bavula2}.
\end{proof}

It is shown in \cite[Theorem 6.50]{Bahturin} that the enveloping algebra
of any Lie algebra (possibly infinite-dimensional) over a field of
characteristic zero does not satisfy any polynomial identity. Similarly,
we have the following.

\begin{theorem}\label{Coulomb-not-PI}
The spherical Coulomb branch algebras $\mathcal{A}(G,N)$ are not PI
algebras.
\end{theorem}

\begin{proof}
This follows from Theorem~\ref{Coulomb-is-Galois-ring},
Lemma~\ref{fundamental-lattices}, and Theorem~\ref{usually-not-PI}.
\end{proof}







\subsection{The Gelfand--Kirillov conjecture for spherical Coulomb branch algebras}

We first record a lemma and then prove that the Gelfand--Kirillov
conjecture holds for spherical Coulomb branch algebras.

\begin{lemma}\label{lemma-Coulomb-1}
The algebra $\mathcal{A}(T,0)$ is isomorphic to
\[
\mathbb{C}[F]\otimes \mathcal{D}(\bar{T}^{\vee}),
\]
and also to a localization of
\[
(\mathbb{C}(\bar{\mathfrak t})*Q)^W,
\]
where $Q$ is the root lattice.
\end{lemma}

\begin{proof}
This follows from a combination of
\cite[Remark 5.23]{BFN} and
\cite[A(i) and Remark A.1]{BFN2}.
\end{proof}

\begin{theorem}\label{Coulomb-GKC}
The Gelfand--Kirillov conjecture holds for all spherical Coulomb branch
algebras $\mathcal{A}(G,N)$. More precisely,
\[
\Frac \mathcal{A}(G,N)
\simeq
\Frac W_n(\mathbb{C}(z_1,\ldots,z_t)),
\]
where
\[
n=\dim\mathfrak t,\qquad t=\dim\mathfrak f .
\]
\end{theorem}

\begin{proof}
By \cite[Remark 5.23]{BFN},  there are isomorphisms between the localizations 
of the algebras
\[
\mathcal{A}(G,N),\qquad
\mathcal{A}(\bar T,N|_{\bar T})^W,
\qquad
\mathcal{A}(\bar T,0)^W .
\]
Hence, using Lemma~\ref{lemma-Coulomb-1}, we obtain
\[
\Frac\mathcal{A}(G,N)
\simeq
\Frac\bigl(\mathbb{C}[F]\otimes
\mathcal{D}(\bar T^\vee)^W\bigr).
\]

Since $T^\vee/W$ is a rational variety, 
\cite[Theorem 1.2]{FS} implies that
\[
\Frac\mathcal{D}(T^\vee)^W
\simeq
\Frac W_n(\mathbb{C}).
\]
The result now follows by adjoining the rational function field
$\mathbb{C}(z_1,\ldots,z_t)$ corresponding to the flavour torus.
\end{proof}

We also obtain the following remarkable consequence.

\begin{corollary}\label{no-CB}
Let $\mathfrak g$ be a finite-dimensional complex simple Lie algebra of
type $B$, $D$, $E$, or $F$. Then its universal enveloping algebra cannot
be realized as a spherical Coulomb branch algebra.
\end{corollary}

\begin{proof}
If it were a spherical Coulomb branch algebra, then by
Theorem~\ref{Coulomb-GKC} its skew field of fractions would be a Weyl
skew field. This contradicts the result of Premet \cite{Premet}.
\end{proof}

The Gelfand--Kirillov conjecture remains open for Lie algebras of type
$C$ and $G_2$, and the Coulomb branch approach provides a possible way
to address this old open problem.

We can also obtain a Poisson analogue of Theorem~\ref{Coulomb-GKC}.
See \cite{Schwarz} for discussions of the Gelfand--Kirillov conjecture
and Noether's problem for Poisson fields.

It was shown in \cite{BFN} that $\mathcal{M}_C(G,N)$ admits a complete
integrable system
\[
\mathcal{M}_C(G,N)\longrightarrow \bar{\mathfrak t}/W
\]
whose generic fibers are isomorphic to $T^\vee$, and that
$\mathcal{M}_C(G,N)$ is birationally equivalent to the cotangent bundle
of $T^\vee/W$. By \cite[Theorem 3.14]{Schwarz}, the Poisson function
field of $T^\vee/W$ is isomorphic to the Poisson field of
$T^*\mathbb{C}^n$; equivalently, it is a standard Poisson field.

More generally, allowing a nontrivial flavour torus $F$ of dimension
$t$, we obtain the following.

\begin{theorem}\label{Coulomb-Poisson}
The Coulomb branch $\mathcal{M}_C(G,N)$ has Poisson function field
\[
\mathbb{C}(x_1,\ldots,x_n,y_1,\ldots,y_n,z_1,\ldots,z_t),
\]
with
\[
\{x_i,y_i\}=1
\]
and all other Poisson brackets equal to zero.
In particular, it is a purely transcendental Poisson field extension
of transcendence degree $2n+t$.
\end{theorem}

\subsection{Ring-theoretic dimensions}

We begin by showing that the algebras $\mathcal{A}(G,N)$ are
LD-stable and by computing their Gelfand--Kirillov dimension.

\begin{theorem}\label{CB-LD-GK}
Let $\mathcal{A}(G,N)$ be a spherical Coulomb branch algebra. Then
$\mathcal{A}(G,N)$ is LD-stable and
\[
\GK\,\mathcal{A}(G,N)=2n+t,
\]
where $n$ is the rank of $G$ and $t$ is the dimension of the flavour
torus.
\end{theorem}

\begin{proof}
It follows from \cite{Zhang2} and \cite{FSS} that Weyl algebras and
their invariant subalgebras are LD-stable. By Theorem
\ref{Coulomb-GKC} and the invariance of the lower transcendence degree
under localization \cite[Proposition 2.1]{Zhang2}, we have
\[
\LD\,\mathcal{A}(G,N)
=
\LD\,\Frac\mathcal{A}(G,N).
\]
The classical computation of the Gelfand--Kirillov transcendence degree
of Weyl fields \cite{Gelfand} therefore gives
\[
\LD\,\mathcal{A}(G,N)=2n+t.
\]

On the other hand, as observed in \cite{lots}, the spherical Coulomb
branch algebra admits a filtration whose associated graded algebra is
the algebra of regular functions on the flavour deformation of the
Coulomb branch $\mathcal{M}_C(G,N)$. This associated graded algebra is a
commutative affine domain. Hence its Gelfand--Kirillov dimension equals
the transcendence degree of its field of fractions. By
Theorem~\ref{Coulomb-Poisson}, this transcendence degree is $2n+t$.
Therefore, by the filtered-graded result
\cite[Corollary to Theorem 1.3]{MStafford},
\[
\GK\,\mathcal{A}(G,N)=2n+t.
\]

Consequently,
\[
\LD\,\mathcal{A}(G,N)
\leq
\Tdeg\,\mathcal{A}(G,N)
\leq
\GK\,\mathcal{A}(G,N)=2n+t,
\]
while
\[
\LD\,\mathcal{A}(G,N)=2n+t.
\]
Hence all these dimensions coincide, and $\mathcal{A}(G,N)$ is
LD-stable.
\end{proof}

\subsection{Gabriel--Rentschler Krull dimension}

We now compute the Krull dimension in the sense of Gabriel--Rentschler.

\begin{theorem}\label{Krull-Coulomb}
Let $\mathcal{A}(G,N)$ be a spherical Coulomb branch algebra with flavour
torus $F$, and set
\[
n=\operatorname{rank}G,\qquad t=\dim F.
\]
Then
\[
\KK(\mathcal{A}(G,N))=n+t .
\]
\end{theorem}

\begin{proof}
By Lemma~\ref{lemma-Coulomb-1}, the algebra
$\mathcal{A}(G,N)$ is a localization of
\[
(\mathbb{C}[\bar{\mathfrak t}]*Q)^W,
\]
where $Q$ is the root lattice acting by translations. Denote this
localization by $\bar{\mathcal D}$.

After a change of lattice basis in
$\operatorname{GL}_{n+t}(\mathbb Z)$, we may write
\[
(\mathbb C[\bar{\mathfrak t}]*Q)^W
\simeq
(\mathbb C[\bar{\mathfrak t}]*\mathbb Z^n)^W,
\]
where $\mathbb Z^n$ acts by translations.

By the results of \cite[Section 7.3]{FO}, this fixed ring can be
identified with
\[
\mathbb C[\mathfrak f]\otimes \mathcal D(T)^W .
\]
The Weyl group $W$ acts by outer automorphisms on the simple ring
$\mathcal D(T)$. Therefore, by Lemma~\ref{krull-invariants} and
\cite[15.3.7]{McConnell},
\[
\KK(\mathcal D(T)^W)=\KK(\mathcal D(T))=n .
\]

Since $\mathbb C[\mathfrak f]$ is a polynomial algebra of dimension $t$,
\cite[6.6.2]{McConnell} gives
\[
\KK((\mathbb C[\bar{\mathfrak t}]*\mathbb Z^n)^W)
=
\KK(\mathbb C[\mathfrak f]\otimes\mathcal D(T)^W)
=
n+t .
\]

By the localization property of Gabriel--Rentschler Krull dimension
(Theorem~\ref{Main-Krull}(i)),
\[
\KK(\mathcal A(G,N))
=
\KK(\bar{\mathcal D})
\leq n+t .
\]

On the other hand, localizing $\bar{\mathcal D}$ further at a suitable
multiplicative set gives
\[
\mathcal D(\mathbb C(x_1,\ldots,x_n))
\otimes
\mathbb C[\mathfrak f].
\]
Since $T/W$ is rational, this algebra has Krull dimension
$n+t$ (cf. \cite[6.6.2, 15.3.10]{McConnell}). Hence another application
of Theorem~\ref{Main-Krull}(i) gives
\[
\KK(\bar{\mathcal D})\geq n+t .
\]

Consequently,
\[
\KK(\mathcal A(G,N))
=
\KK(\bar{\mathcal D})
=
n+t .
\]
\end{proof}

In particular:

\begin{corollary}\label{Krull-W}
 Let $W(\pi)$ be a finite W-algebra of type A (cf. \cite{FMO}). Then its Krull dimension is the (commutative) Krull dimension of its Gelfand-Tsetlin subalgebra, which is $np_1+(n-1)p_2+ \ldots + p_n,$ where $\pi=(p_1,\ldots,p_n)$  
\end{corollary}
\begin{proof}
    It follows from the previous theorem and \cite[Corollary 2.8, Theorem 4.3a)]{WWY} that the Krull dimension is the Krull dimension of the Gelfand-Tsetlin subalgebra, which is a polynomial algebra in $np_1+(n-1)p_2+ \ldots + p_n$ indeterminates (cf. \cite[Section 2]{FMO}). 
\end{proof}

As we saw in the proof of the last Corollary, finite $W$-algebras of type $A$ are spherical Coulomb branch algebras and hence, in particular, all of the results in this section apply to them.

\section*{Acknowledgments}

J. S. would like to acknowledge the invaluable (though intangible) support of J. F. Codas during part of the period in which this manuscript was written.

J.T.H. was supported by Army Research Office grant W911NF-24-1-0058.


\begin{thebibliography}{MMM}

\bibitem[AD06]{AD}
{\sc J. Alev, and F. Dumas},
{\em Opérateurs différentiels invariants et problème de Noether}.
In: Bernstein, J., Hinich, V., Melnikov, A. (eds) Studies in Lie Theory. Progress in Mathematics, vol. 243, Birkhäuser Boston, Boston, MA, 2006, pp. 21--50.





\bibitem[AOV96]{AOvdB}
{\sc J. Alev, A. Ooms and M. Van den Bergh},
{\em A class of counterexamples to the Gelfand-Kirillov conjecture},
Trans. Amer. Math. Soc. {\bf 348} (1996) 1709--1716.






\bibitem[B21]{Bahturin}
{\sc Yu.~A. Bahturin}, {\em Identical Relations in Lie Algebras}. 2nd Ed., De Gruyter Expositions in Mathematics vol. 68, De Gruyter, 2021.


\bibitem[B92]{Bavula}
{\sc V. Bavula},
{\em Generalized Weyl algebras and their representations},
Algebra i Analiz {\bf 4} (1992)
75--97. English translation: St. Petersburg Math. J. {\bf 4} (1993) 71--92.

\bibitem[B95]{Bavula2} {\sc V. Bavula},
{\em Each Schurian algebra is tensor-simple}. Comm. Algebra 23 (1995), no. 4, 1363–1367.






\bibitem[BV00]{BavulaK}
{\sc V. Bavula and F. van Oystaeyen},
{\em Simple holonomic modules over the second Weyl algebra $A_2$},
Adv. Math. 150(1) (2000) 80--116.

\bibitem[BV04]{BO}
{\sc V. Bavula and F. van Oystaeyen},
{\em Simple modules of the Witten-Woronowicz algebra},
J. Algebra {\bf 271} (2004) 827--845.












\bibitem[B81]{Block}
{\sc R. E. Block},
{\em The irreducible representations of the Lie algebra $\mathfrak{sl}(2)$ and of the Weyl algebra},
Adv. Math. {\bf 39}(1) (1981) 69--110.





\bibitem[B75]{Bourbaki2}
{\sc N. Bourbaki,}
{\em Éléments de mathématiques: Groupes et algèbres de Lie}, Ch. VII--VIII, Herman, Paris, 1975.



\bibitem[CH98]{C-B-H}
{\sc W. Crawley-Boevey and M. P. Holland},
{\em Noncommutative deformations of Kleinian singularities},
Duke Math. J. {\bf 92}(3) (1998) 605--635.


\bibitem[BFN18]{BFN} {\sc A. Braverman, M. Finkelberg, H. Nakajima}, {\em Towards a mathematical definition of Coulomb
branches of 3-dimensional N=4 gauge theories, II}, Adv. Theor. Math. Phys. 22 (2018), no. 5,
1071-1147.

\bibitem[BFN19]{BFN2} {\sc A. Braverman, M. Finklberg, H. Nakajima}, {\em Coulomb branches of 3d N=4
 quiver gauge theories and slices in the affine Grassmannian},
Adv. Theor. Math. Phys. 23, No. 1, 75-166 (2019).

\bibitem[BG02]{Brown}
{\sc K.  Brown and K. R. Goodearl},
{\em Lectures on Algebraic Quantum Groups}, Advanced Courses in Mathematics, CRM Barcelona, Birkhäuser Verlag, Basel, 2002.






\bibitem[DF04]{Drensky}
{\sc V. Drensky and E. Formanek},
{\em Polynomial identity rings},
Adv. Courses Math. CRM Barcelona, Birkhäuser Verlag, Basel, 2004.

\bibitem[DOF91]{DFO0}
{\sc Yu. A. Drozd, S. A. Ovsienko, and V. M. Futorny},
{\em On Gel'fand-Zetlin modules},
Proceedings of the Winter School on Geometry and Physics (Srní, 1990), Rend. Circ. Mat. Palermo (2) Suppl. {\bf 26} (1991) 143--147.

\bibitem[DFO94]{DFO}
{\sc Yu, A, Drozd, V. M. Futorny, and S. A. Ovsienko},
{\em Harish-Chandra Subalgebras and Gel'fand-Zetlin Modules}. In: Finite-dimensional Algebras and Related Topics (Ottawa, ON, 1992), 79--93.
NATO Adv. Sci. Inst. Ser. C: Math. Phys. Sci., vol. 424,
Kluwer Academic Publishers Group, Dordrecht, 1994.



\bibitem[EFOS17]{EFOS}
{\sc F. Eshmatov, V. Futorny, S. Ovsienko, and J. Schwarz},
{\em Noncommutative Noether’s problem
for complex reflection groups},
Proc. Amer. Math. Soc. {\bf 145}(12) (2017) 5043--5052.

\bibitem[EG02]{EG}
{\sc P. Etingof and V. Ginzburg},
{\em Symplectic reflection algebras, Calogero-Moser space, and deformed Harish-Chandra homomorphism},
Invent. Math. {\bf 147}(2) (2002) 243--348.

\bibitem[F72]{Faith}
{\sc C. Faith},
{\em Galois subrings of Ore domains are Ore domains}, Bull AMS, {\bf 78}(6) (1972) 1077--1080.

\bibitem[F24]{Fillmore}
{\sc D. Fillmore},
{\em On the Category of Harish-Chandra Block Modules}, 
J. Algebra {\bf 658} (2024) 888--924.
\bibitem[FH25]{FillHart}
{\sc D. Fillmore and J. T. Hartwig},
{\em Subcategories of Module Categories via Restricted Yoneda Embeddings},
arXiv preprint, arXiv:2507.12778[math.RT].


\bibitem[FM78]{FisherM}
{\sc J. W. Fisher, S. Montgomery},
{\em Semiprime skew group rings}, J. Algebra 52 (1978), no. 1, 241--247.



\bibitem[FGRZ20]{FGRZ}
{\sc V. Futorny, D. Grantcharov,  L. E.  Ramirez,  P. Zadunaisky},
{\em Gelfand-Tsetlin theory for rational Galois algebras},
Isr. J. Math. (2020). 

\bibitem[FH14]{FH}
{\sc V. Futorny and J. T. Harwig},
{\em Solution of a $ q$-difference Noether problem and the quantum Gelfand-Kirillov conjecture for $\mathfrak{gl}N$ },
Math. Z. {\bf 276}(no. 1-2) (2014) 1-37. 


\bibitem[FH19]{FH2}
{\sc V. Futorny and J. T. Hartwig},
{\em De Concini–Kac filtration and Gelfand-Tsetlin generators for quantum $gl_N$}, Linear Algebra Appl.568(2019), 173--188.

\bibitem[FKS26]{FKS}
{\sc V. Futorny, P. Koshlukov, J. Schwarz}, {\em Graded identities of the first Weyl algebra and its generalizations}, arXiv preprint, arXiv:2606.11497[Math.RA].


\bibitem[FMO10]{FMO}
{\sc V. Futorny, A. Molev, and S. Ovsienko},
{\em The Gelfand-Kirillov conjecture and Gelfand-Tsetlin modules for finite W-algebras},
Adv. Math. {\bf 223}(3) (2010) 773--796.

\bibitem[FO10]{FO}
{\sc V. Futorny and S. Ovsienko},
{\em Galois orders in skew monoid rings},
J. of Algebra {\bf 324} (2010) 598--630.

\bibitem[FO14]{FO2}
{\sc V. Futorny and S. Ovsienko},
{\em Fibers of characters in Gelfand-Tsetlin categories},
Trans. Amer. Math. Soc. {\bf 366}(8) (2014) 4173--4208.


\bibitem[FS17]{FS0}
{\sc V. Futorny and J. Schwarz},
{\em Galois orders of symmetric differential operators}, 
Algebra Discrete Math. {\bf 23}(1) (2017) 35--46.

\bibitem[FS19]{FS2}
{\sc V. Futorny and J. Schwarz},
{\em Quantum linear Galois orders},
Comm. Algebra {\bf 47}(12) (2019) 5361--5369.


\bibitem[FS20a]{FS}
{\sc V. Futorny and J. Schwarz},
{\em Noncommutative Noether's problem vs classic Noether's problem},
Math. Z. {\bf 295}(3--4) (2020) 1323--1335.

\bibitem[FS20b]{FS3}
{\sc V. Futorny and J. Schwarz},
{\em Algebras of invariant differential operators},
J. Algebra {\bf 542} (2020) 215--229.


\bibitem[FSS21]{FSS}
{\sc V. Futorny, J. Schwarz, and I. Shestakov},
{\em $LD$-stability for Goldie rings},
J. Pure Appl. Algebra {\bf 225}(11) (2021) Paper No. 106741, 14 pp.

\bibitem[GR67]{GR}
{\sc P. Gabriel and R. Rentschler},
{\em Sur la dimension des anneaux et ensembles ordonnés}, C. R. Acad. Sci. Paris Sér. A-B265(1967), A712--A715.





\bibitem[G23]{Gaddis}
{\sc J. Gaddis},
{\em The Weyl algebra and its friends: a survey}, arXiv:2305.01609[math.RA].


\bibitem[GK66]{Gelfand}
{\sc I. M. Gelfand and A. A. Kirillov},
{\em Sur les corps liés aux algèbres enveloppantes des algèbres de Lie},
Inst. Hautes Études Sci. Publ. Math. {\bf 31} (1966) 5--19.

\bibitem[GZ05]{GZ}
{\sc A. Giambruno and M. Zaicev},
{\em Polynomial identities and asymptotic methods}.
Math. Surveys Monogr., vol. 122,
American Mathematical Society, Providence, RI, 2005.

\bibitem[G18]{G18}
{\sc V. Ginzburg},
{\em Nil-Hecke algebras and Whittaker D-Modules}. In: Lie Groups, Geometry, and Representation Theory. Springer, 2018, pp. 137--184.

\bibitem[GKV97]{GKV}
{\sc V. Ginzburg, M. Kapranov, E. Vasserot},
{\em Residue construction of Hecke algebras},
Adv. Math. 128 No 1 (1997) 1--19.


\bibitem[GW04]{GW}
{\sc K. R. Goodearl and R. B. Warfield, Jr.},
{\em An introduction to noncommutative Noetherian rings},
London Math. Soc. Stud. Texts, vol. 61, Cambridge University Press, Cambridge, 2004.



\bibitem[H17]{H0}
{\sc J. T. Hartwig},
{\em The $q$-difference Noether problem for complex reflection groups and quantum OGZ algebras},
Comm. Algebra {\bf 45}(3) (2017) 1166--1176.

\bibitem[H20]{Hartwig}
{\sc J. T. Hartwig},
{\em Principal Galois orders and Gelfand-Zeitlin modules},
Adv. Math. {\bf 359} (2020), 106806, 23pp.

\bibitem[H23]{Hartwig2}
{\sc J. T. Hartwig},
{\em Harish-Chandra modules over Hopf Galois orders},
Int. Math. Res. Not. IMRN {\bf 21} (2023) 18273--18301.

\bibitem[H24]{Hartwig3}
{\sc J. T. Hartwig},
{\em Galois order realization of noncommutative type D Kleinian singularities}, arXiv preprint.\\
\href{https://arxiv.org/abs/2406.20012}{arXiv:2406.20012 [math.RT]}.

\bibitem[H93]{Hodges}
{\sc T. J. Hodges},
{\em Noncommutative Deformations of Type-A Kleinian Singularities},
J. Algebra {\bf 161}(2) (1993) 271--290.

\bibitem[H10]{Hoshi}
{\sc A. Hoshi},
{\em Noether's problem and rationality problem for multiplicative invariant fields: a survey}, arXiv preprint.\\
\href{https://arxiv.org/abs/2010.01517}{arXiv:2010.01517 [math.AG]}.






\bibitem[J21]{Jauch}
{\sc E. C. Jauch,}
{\em An extension of $U(\mathfrak{gl}_n)$ related to the alternating group and Galois orders},
J. Algebra {\bf 569} (2021) 568--594.

\bibitem[J22]{Jauch2}
{\sc E. C. Jauch,}
{\em Maps between standard and principal flag orders},
arXiv preprint.\\
\href{https://arxiv.org/abs/2208.13117}{arXiv:2208.13117 [math.RT]}.



\bibitem[K22]{Joel}
{\sc J. Kamnitzer},
{\em Symplectic resolutions, symplectic duality, and Coulomb branches},
Bull. Lond. Math. Soc. 54, No. 5, 1515-1551 (2022).




\bibitem[KWWY24]{lots}
{\sc J. Kamnitzer, B. Webster, A. Weekes, O. Yacobi},
{\em Lie algebra actions on module categories for truncated shifted Yangians},
Forum Math. Sigma {\bf 12} (2024), Paper No. e18, 69 pp.





\bibitem[KN18]{KN} R. Kodera and H. Nakajima, Quantized Coulomb branches of Jordan quiver gauge
theories and cyclotomic rational Cherednik algebras, String-Math 2016, Proc. Sympos. Pure
Math., vol. 98, Amer. Math. Soc., Providence, RI, 2018, pp. 49–78.

\bibitem[K70]{Krause}
{\sc G. Krause},
{\em On the Krull-dimension of left noetherian left Matlis-rings}, Math.
Zeitschrift 118 (1970), 207--214.

\bibitem[KL00]{KL}
{\sc G. R. Krause and T. H. Lenegan},
{\em Growth of Algebras and Gelfand-Kirillov Dimension},
Grad. Stud. Math., vol. 22, American Mathematical Society, Providence, RI, 2000.





\bibitem[LW23]{LW}
{\sc E. LePage and B. Webster},
{\em Rational Cherednik algebras of $G(\ell, p, n)$ from the Coulomb perspective},
Adv. Math. {\bf 433} (2023) Paper No. 109295, 49 pp.


\bibitem[LMO88]{LMO}
{\sc A. Leroy, J. Matczuk and J. Okni\'nski},
{\em On the Gel'fand-Kirillov dimension of normal localizations and twisted polynomial rings}. In: Perspectives in ring theory (Antwerp, 1987), 205--214. NATO Adv. Sci. Inst. Ser. C: Math. Phys. Sci., 233. Kluwer Academic Publishers Group, Dordrecht, 1988. 

\bibitem[L02]{Levasseur?}
{\sc T. Levasseur},
{\em Krull dimension of the enveloping algebra of a semisimple Lie algebra},
Proc. Amer. Math. Soc. 130 no.12 (2002), 3519--3523.









\bibitem[MT00]{MT}
{\sc V. Mazorchuk, L. Turowska},
{\em On Gelfand-Zetlin modules over Uq(gl n)},  Czechoslovak Journal of Physics 50, 139–144 (2000).


\bibitem[MV21]{MV}
{\sc V. Mazorchuk, E. Vishnyakova},
{\em Harish-Chandra modules over invariant subalgebras in a skew-group ring}, Asian J. Math. 25 (2021), no. 3, 431--454.

\bibitem[MR01]{McConnell} 
{\sc J. C. McConnell and J. C. Robson}
{\em Noncommutative Noetherian rings},
Grad. Stud. Math., vol. 30,
American Mathematical Society, Providence, RI, 2001.




\bibitem[M80]{Montgomery}
{\sc S. Montgomery},
{\em Fixed rings of finite automorphism groups of associative rings},
Lecture Notes in Math., vol. 818, Springer, Berlin, 1980.

\bibitem[MS89]{MStafford}
{\sc J. C. McConnell and J. T. Stafford},
{\em Gel'fand-Kirillov dimension and associated graded modules},
J. Algebra {\bf 125}(1) (1989) 197--214.




\bibitem[N1915]{NoetherX}
{\sc E. Noether},
{\em Der Endlichkeitssatz der Invarianten endlicher Gruppen}, Math. Ann. Vol.77 (1915), pp. 89-92.

\bibitem[P83]{Park}
{\sc J. K. Park}, {\em Skew group rings with Krull dimension}, Math. J. Okayama Univ. 25 (1983), 75--80.

\bibitem[P89]{Passman2}
{\sc D. S. Passman},
{\em Infinite crossed products},
Pure Appl. Math., vol. 135,
Academic Press, Inc., Boston, MA, 1989.



\bibitem[P17]{Petukhov}
{\sc A. Petukhov},
{\em On the Gelfand-Kirillov conjecture for the $W$-algebras attached to the minimal nilpotent orbits},
J. Algebra 470, 289--299 (2017).











\bibitem[P10]{Premet}
{\sc A. Premet},
{\em Modular Lie algebras and the Gelfand-Kirillov conjecture},
Invent. Math. {\bf 181}(2) (2010) 395--420.

\bibitem[R88a]{Rowen}
{\sc L. H. Rowen},
{\em Ring theory. Vol. I.}
Pure Appl. Math., 127
Academic Press, Inc., Boston, MA, 1988.

\bibitem[R88b]{Rowen2} {\sc L. H. Rowen}, {\em Ring theory}, Vol. II. Pure Appl. Math., 128 Academic Press, Inc., Boston, MA, 1988.




\bibitem[S22]{Schwarz}
{\sc J. Schwarz},
{\em A Poisson Noether's Problem and Poisson rationality}, J. Algebra 606 (2022), 195--208.

\bibitem[S25]{SchwarzPan}
{\sc J. Schwarz},
{\em Generalizations of noncommutative Noether's problem},
Journal of Pure and Applied Algebra volume 229, Issue 2, February 2025, 107896.

\bibitem[S26a]{Schwarz3}
{\sc J. Schwarz}. {Characteristic free Galois rings and generalized Weyl algebras}, arXiv preprint, arXiv:2601.10346 [math.RT].


\bibitem[S26b]{Schwarz2}
{\sc J. Schwarz},
{\em Harish-Chandra modules revivisited}, to appear on IMRN.

\bibitem[T22]{Tikaradze}
{\sc A. Tikaradze},
{\em The noncommutative Noether’s problem is almost equivalent to the classical
Noether’s problem}, Adv. Math. 396 (2022), Paper No. 108161, 4 pp.






\bibitem[W24]{Webster}
{\sc B. Webster},
{\em Gelfand-Tsetlin modules in the Coulomb context}, Ann. Represent. Theory 1, No. 3, 393--437 (2024).

\bibitem[WWY20]{WWY} 
{\sc B. Webster, A. Weekes, and O. Yacobi},
{\em A quantum Mirković-Vybornov isomorphism},
Represent. Theory {\bf 24} (2020) 38--84.

\bibitem[W19]{Weekes}
{\sc A. Weekes},
{\em Generators for Coulomb branches of quiver gauge theories}, arXiv:1903.07734.
\href{http://arxiv.org/abs/1903.07734}{arXiv:1903.07734}.


\bibitem[YZ06]{YZ}
{\sc A. Yekutieli and J. J. Zhang},
{\em Homological transcendence degree},
Proc. Lond. Math. Soc. (3) {\bf 93}(1) (2006) 105--137.


\bibitem[Z96]{Zhang}
{\sc J. J. Zhang},
{\em On Gelfand-Kirillov transcendence degree},
Trans. Am. Math. Soc. {\bf 348}(7) (1996) 2867--2899.

\bibitem[Z98]{Zhang2}
{\sc J. J. Zhang},
{\em On lower transcendence degree},
Adv. Math. {\bf 139}(2) (1998) 157--193.

\bibitem[Z73]{Zhelobenko}
{\sc D. P. Zhelobenko},
{\em Compact Lie Groups and Their Representations}, Nauka, Moscow, 1970; Transl. Math. Monogr., vol. 40,
Amer. Math. Soc., Providence, RI, 1973.















\end{thebibliography}
\end{document}